\documentclass[DIV=14,letterpaper]{scrartcl}

\date{}

\usepackage{tikz}
\usetikzlibrary{calc}
\usetikzlibrary{matrix}
\usepackage{amsmath,amssymb,amsfonts,amsthm}
\usepackage[pdftex]{hyperref}
\usepackage{pdfpages}
\usepackage{graphicx}
\usepackage{subfigure}
\allowdisplaybreaks
\usepackage{algorithm,algpseudocode,float}
\usepackage{mathabx}

\newcommand{\V}[1]{\mbox{\boldmath $ #1 $}}

\def \M{\mathbb{M}}

\newcommand{\bey}{\begin{eqnarray}}
\newcommand{\eey}{\end{eqnarray}}

\newcommand{\beq}{\begin{equation}}
\newcommand{\eeq}{\end{equation}}
\theoremstyle{plain}

\theoremstyle{definition}

\newtheorem{exam}{\hspace{6mm}Example}[section]

\title{A Moving Mesh Finite Element Method for Bernoulli Free Boundary Problems}

\author{
Jinye Shen\thanks{School of Mathematics, Southwestern University of Finance and Economics,
Chengdu, Sichuan, China. {\em jyshen@swufe.edu.cn}},
\;
Heng Dai\thanks{School of Mathematics, Southwestern University of Finance and Economics,
Chengdu, Sichuan, China. {\em dh20230810@163.com}},
\; and\,
Weizhang Huang\thanks{Department of Mathematics, University of Kansas, Lawrence, Kansas, U.S.A.
{\em whuang@ku.edu}}
}

\begin{document}
\vskip 1cm
\maketitle

\begin{abstract}
A moving mesh finite element method is studied for the numerical solution of Bernoulli free boundary problems. The method is based on the pseudo-transient continuation with which a moving boundary problem is constructed and its steady-state solution is taken as the solution of the underlying Bernoulli free boundary problem. The moving boundary problem is solved in a split manner at each time step: the moving boundary is updated with the Euler scheme, the interior mesh points are moved using a moving mesh method, and the corresponding initial-boundary value problem is solved using the linear finite element method. The method can take full advantages of both the pseudo-transient continuation and the moving mesh method. Particularly, it is able to move the mesh, free of tangling, to fit the varying domain for a variety of geometries no matter if they are convex or concave. Moreover, it is convergent towards steady state for a broad class of free boundary problems and initial guesses of the free boundary. Numerical examples for Bernoulli free boundary problems with constant and non-constant Bernoulli conditions and for nonlinear free boundary problems are presented to demonstrate the accuracy and robustness of the method and its ability to deal with various geometries and nonlinearities.
\end{abstract}



\newpage

\section{Introduction}

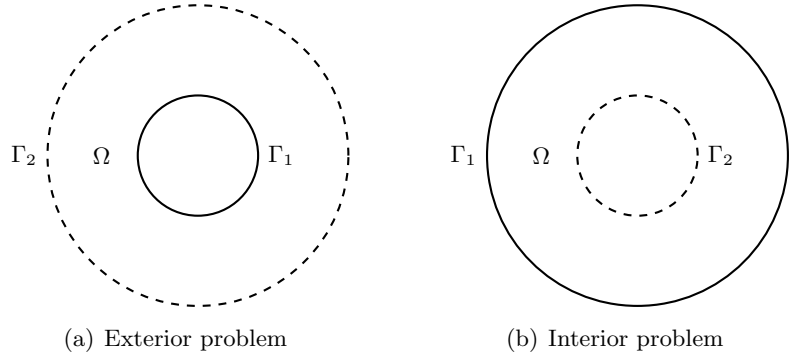
\begin{figure}[tbh]
\centering
\subfigure[Exterior problem]{
\begin{tikzpicture}[scale = 0.8]
\draw[thick,dashed] (0,0)  circle (2.5);
\draw[thick] (0,0)  circle (1);
\node[right] at (1,0) {\footnotesize $\Gamma_1$};
\node at (-1.6,0) {\footnotesize $\Omega$};
\node[left] at (-2.5,0) {\footnotesize $\Gamma_2$};
\end{tikzpicture}
}
\hspace{20pt}
\subfigure[Interior problem]{
\begin{tikzpicture}[scale = 0.8]
\draw[thick] (0,0)  circle (2.5);
\draw[thick,dashed] (0,0)  circle (1);
\node[right] at (1,0) {\footnotesize $\Gamma_2$};
\node at (-1.6,0) {\footnotesize $\Omega$};
\node[left] at (-2.5,0) {\footnotesize $\Gamma_1$};
\end{tikzpicture}
}
\caption{Illustration of the domain for exterior and interior Bernoulli FBPs.}
\label{fig:FBP-1}
\end{figure}

Bernoulli free boundary problems (FBPs) arise in ideal fluid dynamics, optimal insulation, and electro chemistry
\cite{Flucher1997} and serve as a prototype of stationary FBPs.
They have been extensively studied theoretically and numerically; e.g., see
\cite{Alt1981,Brugger2020,Burman2017,Cardaliaguet2002,Crank1984,Eppler2012,Flucher1997,Rabago2020,Weiss2019}.
To be specific, we consider here a typical Bernoulli FBP
\begin{equation}
\label{fbp-1}
\begin{cases}
 - \Delta u  = 0, & \text{ in } \Omega
\\
u  = 1, & \text{ on } \Gamma_1
\\
u  = 0, & \text{ on } \Gamma_2
\\
-\frac{\partial u}{\partial n} =  \lambda, & \text{ on } \Gamma_2
\end{cases}
\end{equation}
where $\Omega$ is a connected domain in $\mathbb{R}^2$ (see Fig.~\ref{fig:FBP-1}),
$\lambda$ is a positive constant, $\Gamma_1 \cup \Gamma_2 = \partial \Omega$, $\Gamma_1$
is given and fixed, and  $\Gamma_2$ is unknown a priori and part of the solution.
We emphasize that the numerical method studied in this work can be applied to more general FBPs
without major modifications, and several such examples are presented in Section~\ref{SEC:numerics-2}.

The Neumann boundary condition in (\ref{fbp-1}) is called the Bernoulli condition.
This condition can be shown to be equivalent to $| \nabla u | = \lambda$ (with the help of the Dirichlet boundary
condition on $\Gamma_2$).
Moreover, the problem is called an exterior (or interior) Bernoulli problem when $\Gamma_2$
is exterior (or interior) to $\Gamma_1$ (cf. Fig.~\ref{fig:FBP-1}).
It is known \cite{Alt1981,Beurling1957,Flucher1997} that an exterior Bernoulli problem has a solution for any $\lambda > 0$
and such a solution is unique and elliptic when the domain enclosed by $\Gamma_1$ is convex.
Loosely speaking, a solution is said to be elliptic (or hyperbolic) if $\Gamma_2$ is getting closer to (or moving away from) $\Gamma_1$
as $\lambda$ increases. On the other hand, an interior Bernoulli problem has a solution only for $\lambda$ large enough and
such solutions are not unique in general. Both elliptic and hyperbolic solutions can co-exist for the same value of $\lambda$
for interior problems.

While the differential equation and boundary conditions are linear, the problem (\ref{fbp-1})
is actually highly nonlinear due to the coupling between $u$ and $\Omega$.
A number of numerical methods have been developed for solving Bernoulli FBPs; e.g., see
a summary of early works for general FBPs \cite[Chapter 8]{Crank1984},
the explicit and implicit Neumann methods \cite{Flucher1997},
a combined level set and boundary element method \cite{Kuster-2007},
shape-optimization-based methods \cite{Eppler2006,Eppler2012,Haslinger2003,Rabago2020},
the cut finite element method \cite{Burman2017},
the quasi-Monte Carlo method \cite{Brugger2020},
the comoving mesh method \cite{SKR2021},
and the singular boundary method \cite{Chen2021}.
A common theme among those methods is trial free boundary and thus iterating between the update of the free boundary
and the solution of the corresponding boundary value problem. Challenges for this approach include
how to choose the initial guess for $\Gamma_2$
to make the iteration convergent and to re-generate or deform the mesh to fit the varying domain.

In this work we shall present a moving mesh finite element method
for the numerical solution of Bernoulli FBPs.
The method is based on the pseudo-transient continuation (e.g., see Fletcher \cite[Section 6.4]{Fletcher1991})
with which we construct an equivalent time-dependent problem (a moving boundary problem or an MBP),
march it until the steady state is reached, and take the steady-state solution
as the solution of Bernoulli FBP (\ref{fbp-1}).
The pseudo-transient continuation is widely used for difficult nonlinear problems in science and engineering
because it can be made convergent for a large class of initial solutions.
Another advantage of using the pseudo-transient continuation is that
the corresponding MBP can be solved readily with boundary-fitted meshes using the finite element method and
the Moving Mesh PDE (MMPDE) method.
The MMPDE method has been developed for general mesh adaptation and movement; e.g., see \cite{HRR94a,HR11}.
It moves the mesh points continuously in time while providing an effective control of mesh quality and concentration.
Most importantly, the method guarantees that the mesh is free of tangling
for any domain (convex or concave) in any spatial dimension \cite{HK2018}. This mesh nonsingularity
is crucial for any mesh-based computation including that for FBPs.
 On the other hand, close attention shluld be paid to the update of the free boundary
where both the gradient of the finite element solution and the normal to the approximate boundary are needed in the computation
of the Bernoulli condition but not defined at boundary vertices in the standard finite element approximation on
a simplicial mesh. Their re-constructions are required and such re-constructions can affect the spatial accuracy of the overall computation.
Two re-construction approaches, (area-)averaging and quadratic least squares fitting, will be discussed.

An outline of this paper is as follows. The pseudo-transient continuation and the corresponding MBP will be described
in Section~\ref{SEC:PTC}. Section~\ref{SEC:MM-FEM} is devoted to the description of the moving mesh FEM.
Numerical examples for Bernoulli FBPs and nonlinear FBPs are presented in Sections~\ref{SEC:numerics} and
\ref{SEC:numerics-2}, respectively. Finally, conclusions and further comments are given
in Section~\ref{SEC:conclusions}.

\section{The pseudo-transient continuation}
\label{SEC:PTC}

For Bernoulli FBP (\ref{fbp-1}), we consider the time-dependent problem
\begin{equation}
\label{fbp-0}
\begin{cases}
 \frac{\partial u}{\partial t} - \Delta u = 0, & \text{ in } \Omega, \text{ for } t > 0
\\
u  = 1, & \text{ on } \Gamma_1
\\
u  = 0, & \text{ on } \Gamma_2
\\
\dot{\Gamma} = - \frac{\partial u}{\partial n} - \lambda, & \text{ on } \Gamma_2 .
\end{cases}
\end{equation}
This system is marched until the steady state is reached and the obtained steady-state solution is taken
as the solution of Bernoulli FBP (\ref{fbp-1}).

\begin{figure}[tbp]
\begin{center}
\begin{tikzpicture}[scale = 0.8]
\draw[thick,dashed] (0,0)  circle (2.5);
\draw[thick] (0,0)  circle (1);
\node[right] at (1,0) {\footnotesize $u = 1$};
\node at (-1.6,0) {\footnotesize $\Omega$};
\node[right] at (2.5,0) {\footnotesize $u = 0$};

\draw[->] (1,2.2913) -- (0.7,1.6039);
\draw[fill=black] (1,2.2913) circle (.5ex);
\node [right, text width = 6cm] at (2.5, 2) {\footnotesize If $(-\frac{\partial u}{\partial n} - \lambda) < 0$,
\\ then $-\frac{\partial u}{\partial n}$ increases and $\dot{\Gamma}$ goes zero.};

\draw[->] (2.3,-0.9798) -- (3.1,-1.3206);
\draw[fill=black] (2.3,-0.9798) circle (.5ex);
\node [right, text width = 6cm] at (2.5,-2) {\footnotesize If $(-\frac{\partial u}{\partial n} - \lambda) > 0$,
\\ then $-\frac{\partial u}{\partial n}$ decreases and $\dot{\Gamma}$ goes zero.};
\end{tikzpicture}
\caption{Illustration of boundary movement for MBP (\ref{fbp-0}).}
\label{fig:boundary-movement}
\end{center}
\end{figure}
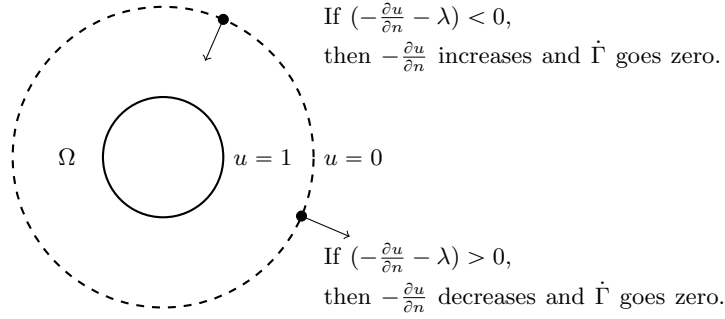

Immediate questions are if MBP (\ref{fbp-0}) has a steady-state solution and whether or not such a solution is stable.
They are related to the asymptotic behavior of solutions of MBPs
and the stability of their steady-state solutions, an area of active research; e.g., see \cite{DuLin2010,FriedmanHu2006,Yamada2020,Zhou2009}.
Unfortunately, none of the available theoretical results seems applicable to (\ref{fbp-0}). Nevertheless, we can gain some
insight from a formal analysis. We take the exterior problem as an example (cf. Fig.~\ref{fig:boundary-movement}).
From the maximum principle,
we know that $u \ge 0$ on $\Omega$ and $\frac{\partial u}{\partial n}|_{\Gamma_2} \le 0$.
Consider a point $\V{x}$ on $\Gamma_2$. If $- \frac{\partial u}{\partial n} - \lambda > 0$ at this point, then $\dot{\Gamma} > 0$ and $\V{x}$ moves outward and farther away from $\Gamma_1$.
Recall that $u|_{\Gamma_1} = 1$ and $u|_{\Gamma_2} = 0$. Thus,
as the distance between $\Gamma_1$ and $\Gamma_2$ increases, $- \frac{\partial u}{\partial n}$ (and therefore,
$\dot{\Gamma} = - \frac{\partial u}{\partial n} - \lambda$) will decrease,
which means that the movement of $\V{x}$ will slow down.
This continues until $- \frac{\partial u}{\partial n} - \lambda$ reaches zero. On the other hand,
if $- \frac{\partial u}{\partial n} - \lambda < 0$, $\V{x}$ will move inward and closer to $\Gamma_1$, which will cause
$- \frac{\partial u}{\partial n}$ to increase and the movement of the boundary point to slow down
until $- \frac{\partial u}{\partial n} - \lambda$ reaches zero. Thus, for either case $\Gamma_2(t)$ will reach a steady state
and so does the domain $\Omega$.
Once $\Omega$ gets close to its steady state, (\ref{fbp-0}) behaves like a parabolic problem with a fixed
domain and its solution will reach steady state too. Thus, (\ref{fbp-0}) reduces to (\ref{fbp-1}).

It is interesting to point out that the use of (\ref{fbp-0}) can also be justified by shape optimization theory.
Indeed, Bernoulli FBPs can be formulated as shape optimization problems; e.g., see \cite{Eppler2006,Eppler2012,Haslinger200301}.
One of the equivalent shape optimization problems for (\ref{fbp-1}) is to minimize the cost function
\begin{equation}
J(\Omega) = \int_\Omega (| \nabla u|^2 + \lambda^2) d \V{x}
\label{fbp-2}
\end{equation}
subject to the PDE constraint
\begin{equation}
\label{fbp-3}
\begin{cases}
 - \Delta u  = 0, & \text{ in } \Omega
\\
u  = 1, & \text{ on } \Gamma_1
\\
u  = 0, & \text{ on } \Gamma_2 .
\end{cases}
\end{equation}
Using shape optimization calculus \cite{Haslinger200301,Sokolowski1992},
we can find the variation of $J(\Omega)$ along a given vector field $\V{V}$ as
\begin{equation}
\delta J(\Omega)[\V{V}] = \int_{\Gamma_2} \left (\lambda^2 - (\frac{\partial u}{\partial n})^2 \right )
\V{V}\cdot \V{n} \; d s .
\label{fbp-4}
\end{equation}
Thus, a descent direction for $J(\Omega)$ is to update $\Gamma_2$ along
\[
- \left (\lambda^2 - (\frac{\partial u}{\partial n})^2 \right )\V{n}
= - (\lambda + \frac{\partial u}{\partial n}) (\lambda - \frac{\partial u}{\partial n}) \V{n} .
\]
The maximum principle implies that the solution to (\ref{fbp-3}) is positive in $\Omega$ and $\frac{\partial u}{\partial n} \le 0$
on $\Gamma_2$. Since $\lambda$ is positive,  we have $(\lambda - \frac{\partial u}{\partial n}) > 0$
and the descent direction is proportional to $- (\lambda + \frac{\partial u}{\partial n}) \V{n}$.
This suggests that $\Gamma_2$ can be updated along $- (\lambda + \frac{\partial u}{\partial n}) \V{n}$, i.e.,
\begin{equation}
\label{Gamma_2_dot}
\dot{\Gamma} = - \frac{\partial u}{\partial n} - \lambda \quad  \text{ on } \Gamma_2.
\end{equation}
This gives the boundary velocity in (\ref{fbp-0}).
Interestingly, trial free boundary methods (e.g., see \cite[Chapter 8]{Crank1984})
can be interpreted as a time discretization of the above equation.
Moreover, Sunayama et al. \cite{SKR2021} use
\begin{equation}
\label{fbp-5}
\begin{cases}
 - \Delta u  = 0, & \text{ in } \Omega
\\
u  = 1, & \text{ on } \Gamma_1
\\
u  = 0, & \text{ on } \Gamma_2
\\
\dot{\Gamma} = - \frac{\partial u}{\partial n} - \lambda, & \text{ on } \Gamma_2.
\end{cases}
\end{equation}
The difference between this system and (\ref{fbp-0}) is that
the heat equation, instead of the Laplace equation, is used in (\ref{fbp-0}).
With (\ref{fbp-0}), we can take full advantages of the pseudo-transient continuation
in the numerical solution. Particularly, we can use automatic time stepsize selection procedures and extend
the developed numerical method to more general FBPs (including nonlinear FBPs) without major modifications
(cf. Section~\ref{SEC:numerics-2}).

\section{A moving mesh finite element solution for MBPs}
\label{SEC:MM-FEM}

In this section, we describe a moving mesh finite element method for solving MBP (\ref{fbp-0}).
The method solves (\ref{fbp-0}) in a splitting manner at each time step: updates the boundary using the Euler scheme,
moves the interior mesh points using the MMPDE moving mesh method, and integrates the underlying
initial-boundary value problem using a Runge-Kutta scheme and linear finite elements. The method has been
used in \cite{NgoH-2019} for solving the porous medium equation. The method can be used for general MBPs although
it is described here only for (\ref{fbp-0}).

\begin{algorithm}[htbp]
\caption{Moving mesh FEM for (\ref{fbp-0})}
\label{MMFEM-euler}
\begin{itemize}
\item[0.] Assume that $\mathcal{T}_h^n$ and $u_h^n$ at $t = t_n$ are known.
\item[1.] \textbf{Boundary update.} Update the mesh vertices on $\Gamma_2$ using the  Euler scheme,
	\begin{equation}\label{Forward-1}
	\V{x}_i^{n+1} = \V{x}_i^{n} + \Delta t_n \left.\left( - \nabla u_h^n \cdot \V{n}  -\lambda \right)\V{n}
	\right |_{\V{x}_i^n}, \quad \forall \V{x}_i^n \in \Gamma_2^n .
	\end{equation}
	Denote by $\Gamma_2^{n+1}$ the updated boundary and by $\tilde{\mathcal{T}}_h^{n+1}$ the mesh with the updated
	boundary. 	Thus, the vertices of $\tilde{\mathcal{T}}_h^{n+1}$ consist of the boundary vertices on $\Gamma_2^{n+1}$ and
	$\Gamma_1$ and the interior vertices of $\mathcal{T}_h^n$.
	 Notice that the Euler update (\ref{Forward-1}) generally will not result in an even distribution of the boundary vertices along
	the boundary. They can be made more evenly distributed in the next step (the mesh movement step) by allowing the boundary vertices
	to slide along the boundary.

\item[2.] \textbf{Movement of interior mesh vertices.} Generate the new mesh $\mathcal{T}_h^{n+1}$ for $\Omega^{n+1}$
	by moving the vertices of $\tilde{\mathcal{T}}_h^{n+1}$ using the MMPDE moving mesh method.
	The detail is given in Subsection~\ref{sec:MMPDE}.

\item[3.] \textbf{Solution of the initial-boundary value problem.}  Solve the IBVP
	\begin{equation}
	\label{mbp-2}
	\begin{cases}
	\frac{\partial u}{\partial t} - \Delta u = 0, & \text{ in } \Omega(t)
	\\
	u  = 1, & \text{ on } \Gamma_1
	\\
	u  = 0, & \text{ on } \Gamma_2(t)
	\end{cases}
	\end{equation}
	on the moving mesh $\mathcal{T}_h(t)$ defined as the linear interpolation between $\mathcal{T}_h^{n}$ and
	$\mathcal{T}_h^{n+1}$, i.e.,
	\begin{equation}
	\V{x}_i(t) = \frac{t_{n+1}-t}{\Delta t_n} \V{x}_i^{n} + \frac{t-t_{n}}{\Delta t_n} \V{x}_i^{n+1}, \quad i = 1, ..., N_v,
	\quad t \in [t_n, t_{n+1}].
	\label{x-0}
	\end{equation}
	 In this step, the domain moves from $\Omega^n$ to $\Omega^{n+1}$ and is considered known
	 (as specified by the meshes $\mathcal{T}_h^n$ and $\mathcal{T}_h^{n+1}$).
	 Piecewise linear finite elements and a fifth-order implicit Runge-Kutta scheme
	are employed for the spatial and temporal discretization of the IBVP, respectively. The detail
	is given in Subsection~\ref{sec:FEM}.
\end{itemize}
\end{algorithm}

\subsection{The overall procedure of the moving mesh FEM}
\label{sec:3.1}

Denote the time instants by $t_n$, $n = 0, 1, ...$ and
the corresponding time steps by $\Delta t_n = t_{n+1} - t_n$.
We assume that the moving domain $\Omega(t)$
is partitioned into/approximated by a moving triangular mesh $\mathcal{T}_h(t)$ that has $N_v$ vertices
(denoted by $\V{x}_i(t), i=1,\ldots,N_v$), $N$ elements, and a fixed connectivity. The domain and mesh at $t_n$
will be denoted by $\Omega^n$ and $\mathcal{T}_h^n$, respectively.
The goal of the moving mesh FEM is to generate a new mesh $\mathcal{T}_h^{n+1}$  and
a new numerical solution $u_h^{n+1}$ at any given time $t=t_{n+1}$. The method
contains three basic steps and its overall procedure is given in Algorithm~\ref{MMFEM-euler}.
Since the boundary movement and the update of the physical solution are split and performed sequentially,
the method is expected to be first-order in time. Moreover, the physical PDE is discretized spatially with linear finite elements
and we expect the method to be second-order in space.
 Notice that the lower-oder convergence in time for the moving mesh FEM is not a concern here
since our goal is to obtain a steady-state solution of IBVP (\ref{fbp-0}) and the accuracy of such a steady-state solution
is determined only by spatial discretization.
 Moreover, the use of a triangular mesh for the moving domain gives a piecewise linear approximation to the moving boundary,
which is sufficiently accurate for a second-order numerical approximation for the underlying FBP. For higher-order accuracy, however,
a higher-order mesh, such as the one with a piecewise quadratic approximation to the boundary, has to be used.

Notice that $\nabla u_h^n$ is used in (\ref{Forward-1}).
Since the FE approximation $u_h^n$ is only piecewise linear, its gradient is not defined at vertices (including boundary vertices).
It can be approximated as an area-weighted average of the gradient on the neighboring elements;
e.g. see Murea and Hentschel \cite{Murea-2005} and Ngo and Huang \cite{NgoH-2019}.
Another technique is least squares fitting.
For example, a quadratic polynomial can be formed by fitting the values of $u_h^n$
at the neighboring vertices and differentiated to obtain an approximate gradient.
Furthermore, recently Sturm \cite{Sturm2015} and Sunayama et al. \cite{SKR2021} proposed to define a mesh velocity field
on the whole domain by solving a Laplace boundary value problem with Dirichlet/Robin boundary conditions.
In our computation we use the quadratic least squares fitting and compare it with the area-weighted averaging technique.
Numerical results show that the quadratic least squares fitting can lead to second-order convergence in space
whereas the area-weighted averaging seems to give only first-order convergence.

The unit outward normal $\V{n}$ to the boundary in (\ref{Forward-1}) is not defined at boundary vertices either.
It can be computed either as the average of the unit outward normals on the edges connecting $\V{x}_i^n$ or through the quadratic
least squares fitting.
Numerical results show that the averaging approach maintains the second-order spatial convergence of the method
and thus this approach is used in our computation.
Generally speaking, we can expect this to work when the boundary is sufficiently smooth and the mesh is sufficiently fine.

The mesh $\tilde{\mathcal{T}}_h^{n+1}$, formed after the update of $\Gamma_2$,
is required to be nonsingular (i.e., free of tangling). This can be achieved when
$\Gamma_2$ is sufficiently smooth, the mesh is sufficiently fine, and $\Delta t_n$ is
sufficiently small; e.g., see the analysis of conforming triangulation for moving domains by
Rangarajan and Lew \cite{Lew-2013, Lew-2014}.
Generally speaking, this nonsingularity requirement
of $\tilde{\mathcal{T}}_h^{n+1}$ places a restriction on the maximum time step allowed in the computation.
To see this, it is reasonable to expect that the mesh $\tilde{\mathcal{T}}_h^{n+1}$ stays nonsingular if the boundary vertices
move no more than $a_h^n/2$ over a step,
where $a_h^n$ is the minimum element height of $\mathcal{T}_h^{n}$.
From (\ref{Forward-1}), we have
\begin{equation}
\label{dt-1}
\Delta t_n \le \frac{a_h^n}{ 2 \left (\lambda + \max\limits_{\V{x}_i^n
\in \Gamma_2^n} \left |\frac{\partial u_h^n}{\partial n} (\V{x}_i^n)\right | \right )} .
\end{equation}
The above inequality implies $\Delta t_n = \mathcal{O}(h)$ if the mesh is close to being uniform and $\frac{\partial u_h^n}{\partial n}$
is bounded. Generally speaking, this is not a serious restriction on the time step.
In practice, the nonsingularity of $\tilde{\mathcal{T}}_h^{n+1}$ is checked at each time step  by computing the minimum
height of the mesh elements that should stay away from zero for any nonsingular mesh; the interested read is referred
to the analysis in  \cite{HK2018}.
When $\tilde{\mathcal{T}}_h^{n+1}$
is found to be singular, $\Delta t_n$ is reduced and the boundary is re-computed. This process is repeated until
$\tilde{\mathcal{T}}_h^{n+1}$ is nonsingular.

In Step 2 of Algorithm~\ref{MMFEM-euler},
the new mesh $\mathcal{T}_h^{n+1}$ is generated from the initial mesh $\tilde{\mathcal{T}}_h^{n+1}$
using the MMPDE method.
It has been proven in \cite{HK2018} that the MMPDE method produces a nonsingular mesh
for any (convex or concave) domain in any spatial dimension
if the initial mesh is nonsingular. Thus, the nonsingularity of $\tilde{\mathcal{T}}_h^{n+1}$ implies
the nonsingularity of $\mathcal{T}_h^{n+1}$. More detail of the MMPDE method is given in Subsection~\ref{sec:MMPDE}.

\subsection{Finite element discretization of PDEs on moving meshes}
\label{sec:FEM}

In this subsection we describe the linear FE solution of the IBVP (\ref{mbp-2}) on the moving mesh
$\mathcal{T}_h(t)$ from $t_n$ to $t_{n+1}$. We use the quasi-Lagrange approach (e.g., see \cite{HR11})
where the mesh is considered to move continuously in time (cf. (\ref{x-0})). The nodal velocities are given by
\begin{equation}
\label{x-2}
\dot{\V{x}}_i(t)=\frac{\V{x}_i^{n+1}-\V{x}_i^n}{t_{n+1}-t_n},\quad i=1,\ldots,N_v, \quad t \in (t_n, t_{n+1}).
\end{equation}

Denote the piecewise linear basis function associated with vertex $\V{x}_i$ by $\phi_i(\V{x},t)$.
It depends on $t$ through the movement of vertices. It is not difficult to show
\begin{equation}
\label{x-3}
\frac{\partial \phi_i}{\partial t}=-\nabla \phi_i\cdot \dot{\V{X}},
\end{equation}
where $\dot{\V{X}}$ is the piecewise linear velocity function defined as
\[
\dot{\V{X}}(\V{x},t)=\sum_{i=1}^{N_v}\dot{\V{x}}_i(t)\phi_i(\V{x},t) .
\]
If we arrange the vertices in such a way that the first $N_{vi}$ vertices are the interior vertices,
we can express the linear finite element spaces as
\begin{align*}
& V_h(t)=\text{span}\{\phi_1(\cdot,t),\ldots,\phi_{N_{v}}(\cdot,t)\} \cap \{ v_h |_{\Gamma_1} = 1, \;  v_h |_{\Gamma_2} = 0\},
\\
& V_h^0(t)=\text{span}\{\phi_1(\cdot,t),\ldots,\phi_{N_{vi}}(\cdot,t)\} .
\end{align*}
 Notice that $V_h(t)$ and $V_h^0(t)$ are subspaces of Sobolev spaces $H^1(\Omega)$ and $H^1_0(\Omega)$, respectively.
Then the linear finite element approximation of (\ref{mbp-2}) is to find $u_h(t)\in V_h(t)$, $t > 0$, such that
\begin{align}
 & \int_{\Omega} \frac{\partial u_h}{\partial t}\psi \,d\V{x} +\int_{\Omega} \nabla \psi\cdot \nabla u_h  \,d\V{x}
 =0,\qquad \forall \psi \in V_h^0(t).
\label{V-F}
\end{align}
Expressing $u_h$ as
\begin{equation}\label{u_h}
u_h(\V{x},t)=\sum_{i=1}^{N_{v}}u_i(t)\phi_i(\V{x},t),
\end{equation}
differentiating it with respect to $t$, and using (\ref{x-3}), we obtain
\[
\frac{\partial u_h}{\partial t}=\sum_{i=1}^{N_{v}}\frac{d u_i}{d t}\phi_i(\V{x},t)+\sum_{i=1}^{N_{v}}u_i(t)\frac{\partial \phi_i}{\partial t}
= \sum_{i=1}^{N_{v}}\frac{d u_i}{d t}\phi_i(\V{x},t)-\nabla u_h\cdot \dot{\V{X}} .
\]
Substituting the above equation into \eqref{V-F} and taking $\psi=\phi_j,\;j=1,\ldots,N_{vi}$ successively, we get
\begin{align}
  &\sum_{i=1}^{N_{v}}\left(\int_{\Omega}\phi_i \phi_j \,d\V{x}\right)\frac{d u_i}{d t}
  - \int_{\Omega} \nabla u_h\cdot \dot{\V{X}} \phi_j \, d\V{x}
  +\int_{\Omega} \nabla \phi_j\cdot \nabla u_h  \,d\V{x} = 0, \qquad j =1,\ldots,N_{vi}.
 \label{V-F-1}
\end{align}
This system, together with the boundary conditions, can be cast into a matrix form as
\begin{equation}
\label{Matrix}
  \V{B}(\V{X})\dot{\V{U}}=F(\V{U},\V{X},\dot{\V{X}}) ,
\end{equation}
where $\V{U} = (u_1, ..., u_{N_v})^T$ and $\V{X} = (\V{x}_1, ..., \V{x}_{N_v})^T$.
In principle, any time marching scheme can be used to integrate the above system of ordinary differential equations.
We use the fifth-order implicit Radau IIA Runge-Kutta scheme with variable time step. The selection of time step
is based on a two-step error estimator developed by Gonzalez-Pinto et al. \cite{Montijano-2004} and
the relative and absolute tolerances are chosen as $10^{-6}$ and $10^{-8}$, respectively, in our computation.

 We recall that the moving mesh FEM described in Algorithm~\ref{MMFEM-euler} is first-order
in time overall due to its splitting implementation and Euler update of the moving boundary.
As such, it is more consistent to use a first-order scheme for integrating (\ref{Matrix}).
The choice of the fifth-order implicit Radau IIA Runge-Kutta scheme in our computation is mainly based on the convenience:
the scheme and related time step selection have been implemented in MMPDElab \cite{MMPDElab},
a publicly available Matlab package for adaptive mesh movement and finite element computation in one, two, and three
dimensions. MMPDElab was used in our computation for integrating (\ref{Matrix}) and generating moving meshes (see the next subsection).

\subsection{The MMPDE moving mesh method}
\label{sec:MMPDE}

We use the MMPDE moving mesh method to generate the new mesh $\mathcal{T}_h^{n+1}$
for $\Omega^{n+1}$ starting from $\tilde{\mathcal{T}}_h^{n+1}$. The method has been
developed (e.g., see \cite{HK2015,HRR94a,HR11}) for general mesh adaptation and movement.
It uses the so-called moving mesh PDE (or moving mesh equations in discrete form) to move vertices
continuously in time and in an orderly manner in space.
A key idea of the MMPDE method is viewing any nonuniform mesh as a uniform one in some Riemannian metric
specified by a tensor $\M = \M(\V{x},t)$. For our current situation, the solution of (\ref{fbp-0})
is smooth in space and mesh adaptation is not necessary. Moreover, (\ref{dt-1}) suggests that
a uniform mesh may provide an advantage over nonuniform meshes since it allows a larger time step.
For these reasons, we take $\M = \mathbb{I}$ (the identity matrix) and try to make the mesh as uniform
as possible.

It is known (e.g., see \cite{HK2015,HR11}) that a uniform mesh satisfies the following equidistribution and alignment conditions,
\begin{align}
  |K| &=\frac{\sigma_h}{N},\;\qquad \forall K\in \mathcal{T}_h\;\label{C-1}\\
  \frac{1}{2}\text{trace}\left((F'_{K})^{-1} (F'_{K})^{-T}\right)&
  =\det\left((F'_{K})^{-1}(F'_{K})^{-T}\right)^{\frac{1}{2}},\;\qquad \forall K\in \mathcal{T}_h\;\label{C-2}
\end{align}
where $|K|$ is the area of $K$, $F'_K$ is the Jacobian matrix of the affine mapping
$F_K: \hat{K}\to K$, $\hat{K}$ is the reference element taken as an equilateral triangle with unit area, and
$\sigma_h=\sum_{K\in \mathcal{T}_h}|K|$.
The condition \eqref{C-1} requires all elements to have the same size while \eqref{C-2} requires every element $K$
to be similar to $\hat{K}$.
Since $\hat{K}$ is taken as an equilateral triangle, these conditions actually tempt to make all elements as uniform
and equilateral as possible.
An energy function associated with these conditions is given by
\begin{align}\label{I-h}
  I_h=&\frac{1}{3}\sum_{K\in \mathcal{T}_h}|K| \text{trace}\left((F'_{K})^{-1}(F'_{K})^{-T}\right)^{\frac{3}{2}}
  + \frac{2^{\frac{3}{2}}}{3}\sum_{K\in \mathcal{T}_h}|K|\left(\det(F'_K)\right)^{-\frac{3}{2}}.
\end{align}
This function is a Riemann sum of a continuous functional developed based on mesh equidistribution and alignment (e.g., see \cite{HR11}).

The energy function $I_h$ is a function of the coordinates of the vertices of $\mathcal{T}_h$, i.e.,
$I_h = I_h(\V{x}_1, ..., \V{x}_{N_v})$.
An approach for minimizing this function is to integrate the gradient system of $I_h$. Thus, we define the moving mesh equations as
\begin{equation}\label{x-1}
   \frac{d \V{x}_i}{d t}=-\frac{1}{\tau}\frac{\partial I_h}{\partial \V{x}_i},\qquad i=1,\ldots,N_v
 \end{equation}
where $\tau>0$ is a parameter used to adjust the time scale of mesh movement.
The analytical expression of the derivative of $I_h$ with respect to $\V{x}_i$ can be found using scalar-by-matrix
differentiation \cite{HK2015}. Using this expression, we can rewrite (\ref{x-1}) as
\begin{equation}
 \frac{d \V{x}_i}{d t} = \frac{1}{\tau} \sum\limits_{K \in \omega_i} |K| \V{v}_{i_K}^K ,\qquad
 i = 1, ..., N_v
\label{mmpde-1}
\end{equation}
where $\omega_i$ is the element patch associated with vertex $\V{x}_i$ and $\V{v}_{i_K}^K$ is the local mesh velocity
contributed by element $K$ to the vertex $\V{x}_i$.
Define the edge matrices of $K$ and $\hat{K}$ as
$E_K=[\V{x}_1^K-\V{x}_0^K,\V{x}_2^K-\V{x}_0^{K}]$ and
$\hat{E} =[\V{\xi}_1-\V{\xi}_0,\V{\xi}_2-\V{\xi}_0]$, respectively, where $\V{x}_0^K,\; \V{x}_1^K,\; \V{x}_2^K$
and $\V{\xi}_0,\; \V{\xi}_1,\; \V{\xi}_2$ are the coordinates of the vertices of $K$ and $\hat{K}$.
Let $\mathbb{J}=\hat{E} E_K^{-1}$.
Then, the local mesh velocities are given by
\begin{align*}
& \begin{bmatrix}
  \V{v}_1^K,
 \V{v}_2^K
\end{bmatrix}^T =- G E_K^{-1} + E_{K}^{-1}\frac{\partial G}{\partial \mathbb{J}} \hat{E} E_K^{-1}
+ \frac{\partial G}{\partial \det(\mathbb{J})}\frac{\det(\hat{E})}{\det(E_K)} E_K^{-1},
\quad \V{v}_0^{K}=- \left(\V{v}_1^{K}+\V{v}_2^{K}\right),
\\
&  G(\mathbb{J},\det(\mathbb{J}))=\frac{1}{3}
(\text{trace}(\mathbb{J} \mathbb{J}^{T}))^{\frac{3}{2}}
+\frac{2^{\frac{3}{2}}}{3} \left(\det(\mathbb{J})\right)^{\frac{3}{2}},
\\
&  \frac{\partial G}{\partial \mathbb{J}}=(\text{trace}(\mathbb{J}\mathbb{J}^{T}))^{\frac{1}{2}}\mathbb{J}^{T},\\
&  \frac{\partial G}{\partial \det(\mathbb{J})}=2^{\frac{1}{2}} \det(\mathbb{J})^{\frac{1}{2}}.
\end{align*}
The nodal velocity needs to be modified at boundary vertices. For fixed boundary
vertices, $\frac{d \V{x}_i}{d t}$ should be set to be zero. If $\V{x}_i$ is allowed to slide along the boundary,
the component of $\frac{d \V{x}_i}{d t}$ in the normal direction of the boundary should be set to be zero.
 Allowing the boundary vertices to slide along the boundary is useful in making
them more evenly distributed.
In our computation,  the boundary vertices on the fixed boundary $\Gamma_1$ are fixed while those on the moving boundary
$\Gamma_2$ are allowed to slide along the boundary.
Moreover, the Matlab ODE solver \textit{ode15s} (a variable-step, variable-order solver
based on the numerical differentiation formulas of orders 1 to 5) is used for integrating (\ref{mmpde-1}), with
the Jacobian matrix approximated by finite differences.  The MMPDE method has been implemented
in the Matlab package MMPDElab \cite{MMPDElab}.

It is worth pointing out that there exist other adaptive moving mesh methods; e.g., see
textbooks/reviews \cite{Bai94a,Baines-2011,BHR09,HR11,Tang05} and references therein.
The interested reader is also referred to some recent works on moving mesh methods
\cite{BSD-2019,GuYaguang2022,WTT-2017,ZhangHQ-2022}.

\section{Numerical examples of Bernoulli FBPs}
\label{SEC:numerics}

We now present numerical results obtained for four examples of Bernoulli FBPs with the moving mesh FEM described
in the previous section. Unless stated otherwise,
we use $\tau = 10^{-5}$, $\Delta t_{max} = 0.001$, the zero initial condition
$u(\V{x},0) = 0$, and the quadratic least squares fitting approach for computing $\nabla u_h$ needed in boundary update.
The computation is stopped when the ratio of the current maximum boundary velocity with the initial maximum
boundary velocity is below $10^{-4}$.

\begin{figure}[htbp]
\centering
\subfigure[$t=0$]{
\includegraphics[width=5cm]{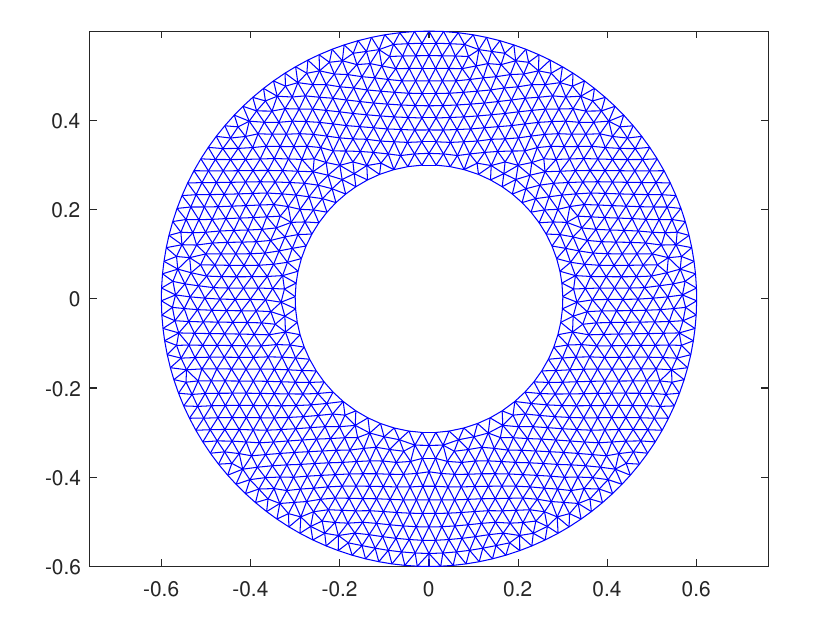}
}
\quad
\subfigure[$t=0.15$]{
\includegraphics[width=5cm]{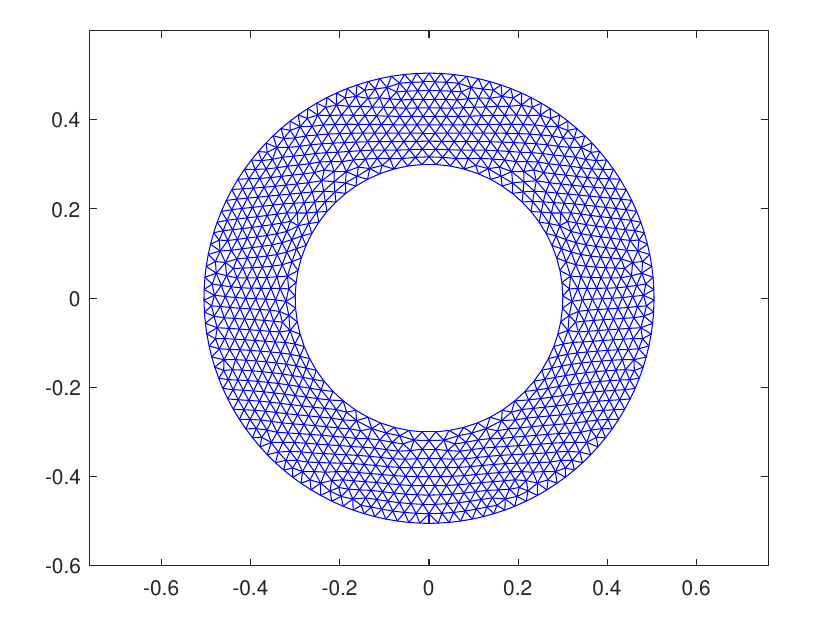}
}
\\
\quad
\subfigure[$t=0.3$]{
\includegraphics[width=5cm]{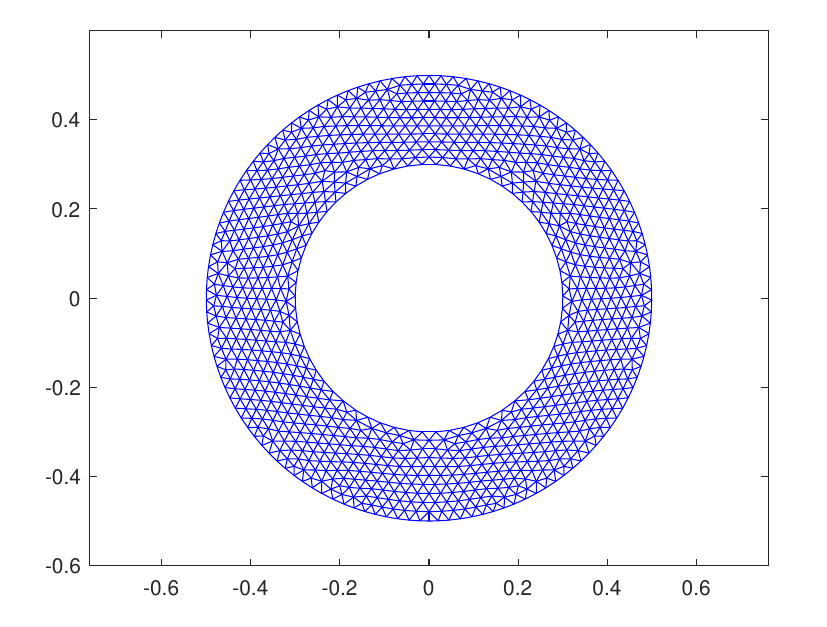}
}
\quad
\subfigure[$t=0.456$]{
\includegraphics[width=5cm]{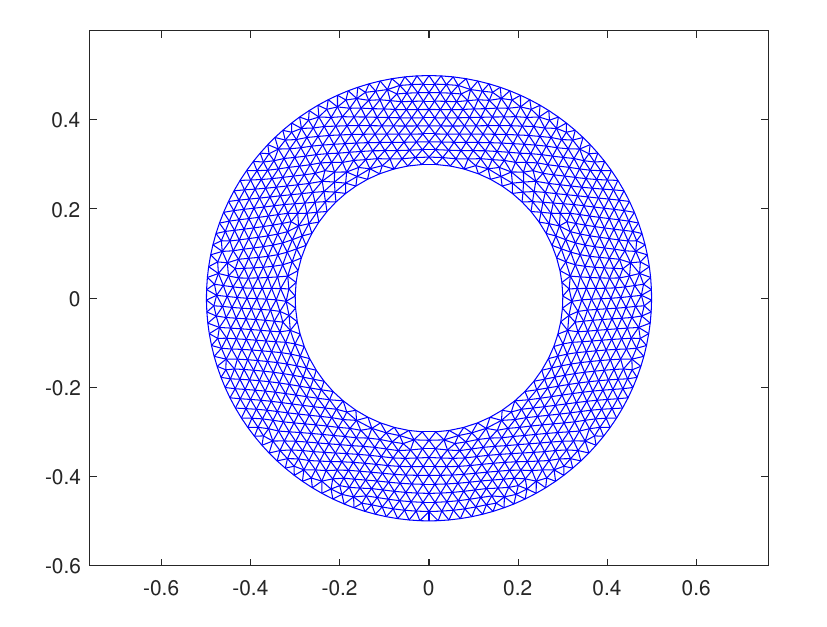}
}
\caption{Example~\ref{fbp-ex1}. The mesh of $N = 1998$ is plotted at $t=0$, 0.15, 0.3, and 0.456 for $\lambda = -2/\ln(0.6)$.}
\label{fig:fbp-ex1-2}
\end{figure}

\begin{figure}[ht]
\centering
\includegraphics[width=4.7cm]{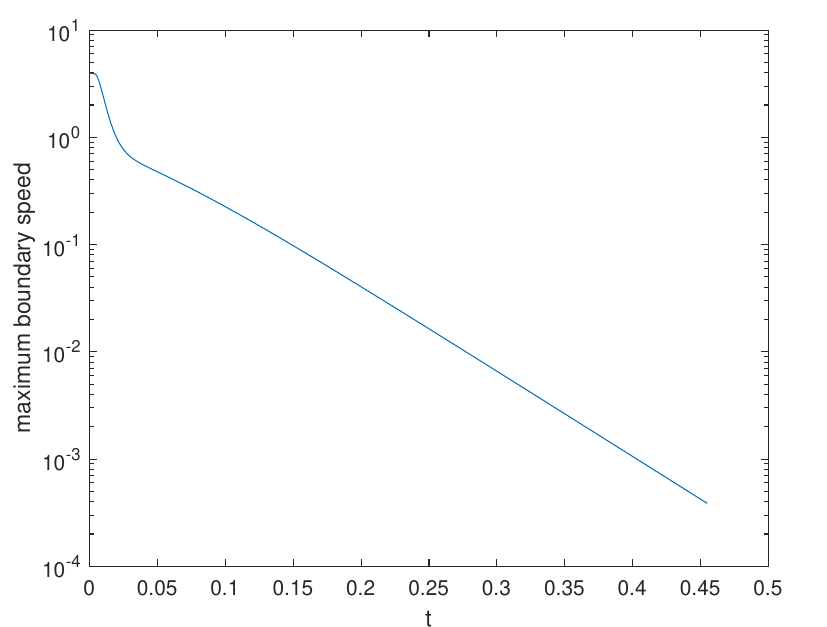}
\caption{Example~\ref{fbp-ex1}. The maximum boundary velocity is plotted as a function of time for $\lambda = -2/\ln(0.6)$ and $N =1998$.}
\label{fig:fbp-ex1-3}
\end{figure}

\begin{figure}[h]
\centering
\includegraphics[width=4.7cm]{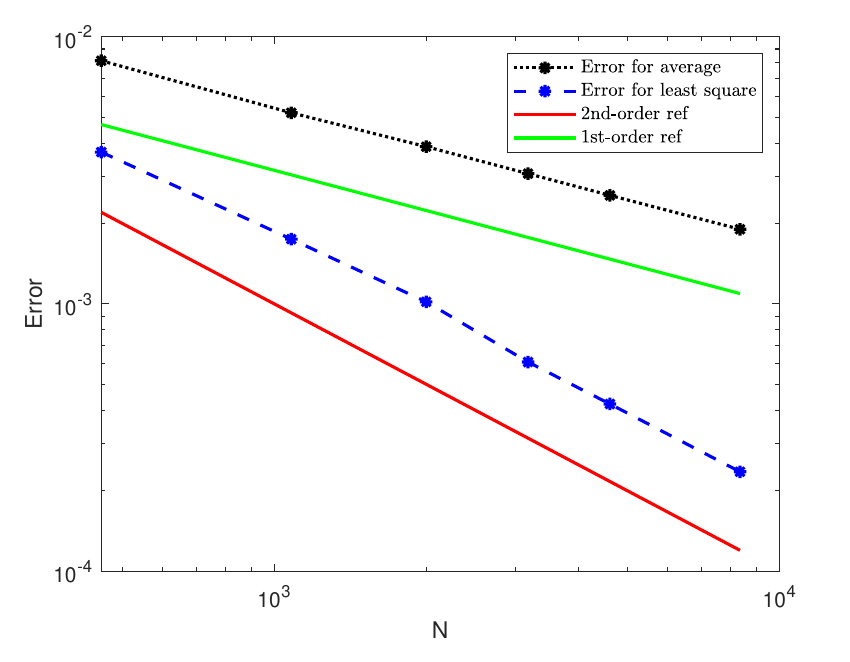}
\caption{Example~\ref{fbp-ex1}. The error in the boundary location is plotted as a function of $N$
(the number of elements in the mesh) for two strategies (the quadratic least squares fitting and area-weighted averaging)
for computing solution gradient used in boundary update.}
\label{fig:fbp-ex1-1}
\end{figure}

\begin{figure}[htbp]
\centering
\subfigure[$t=0$]{
\includegraphics[width=4.3cm]{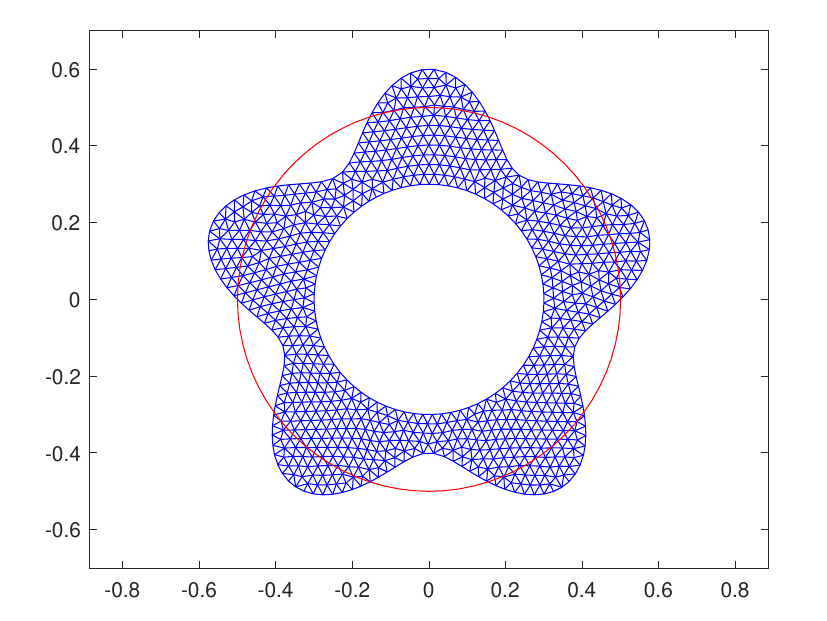}
}
\quad
\subfigure[$t=0.01$]{
\includegraphics[width=4.3cm]{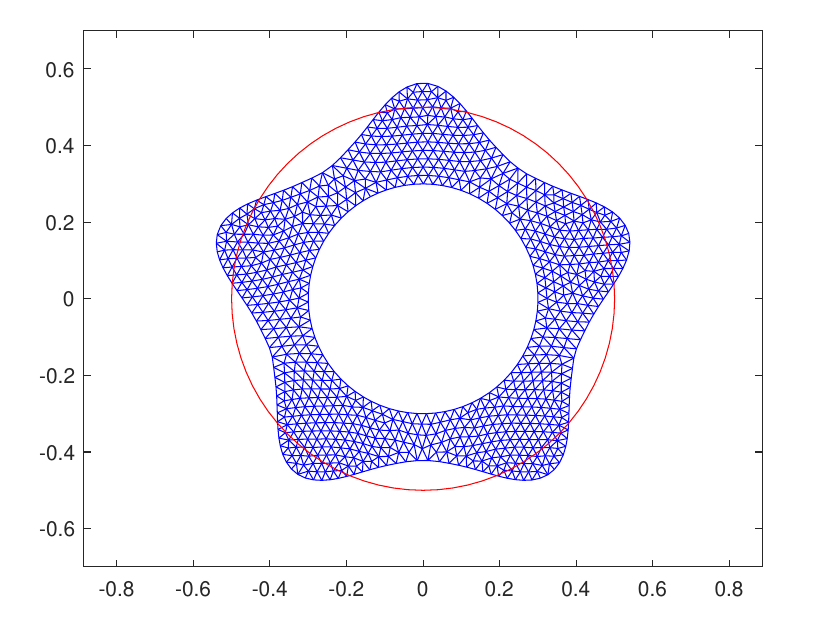}
}
\quad
\subfigure[$t=0.02$]{
\includegraphics[width=4.3cm]{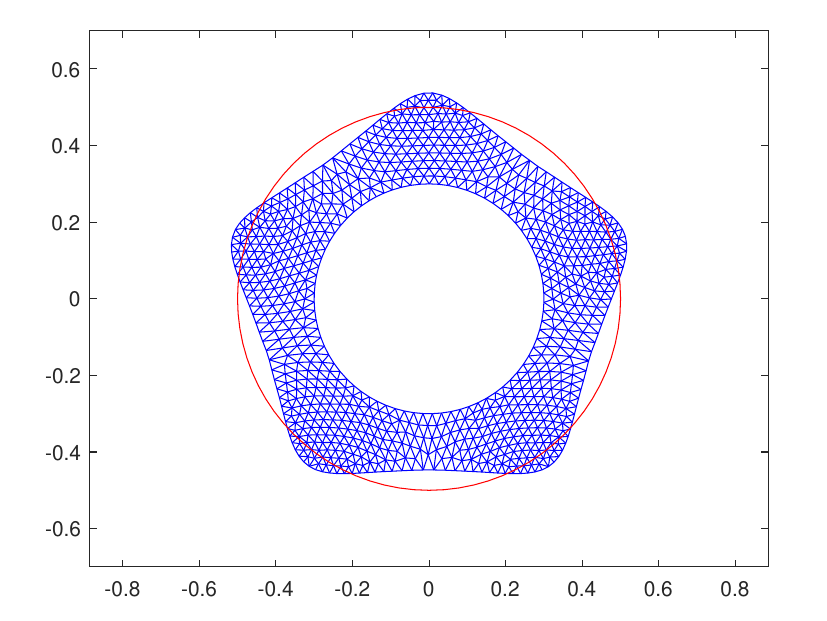}
}
\quad
\subfigure[$t=0.05$]{
\includegraphics[width=4.3cm]{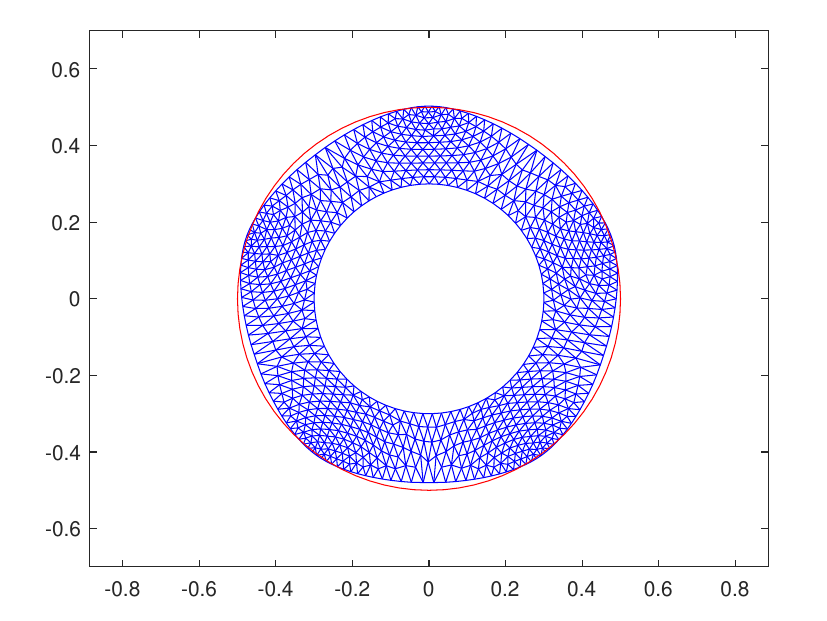}
}
\quad
\subfigure[$t=0.1$]{
\includegraphics[width=4.3cm]{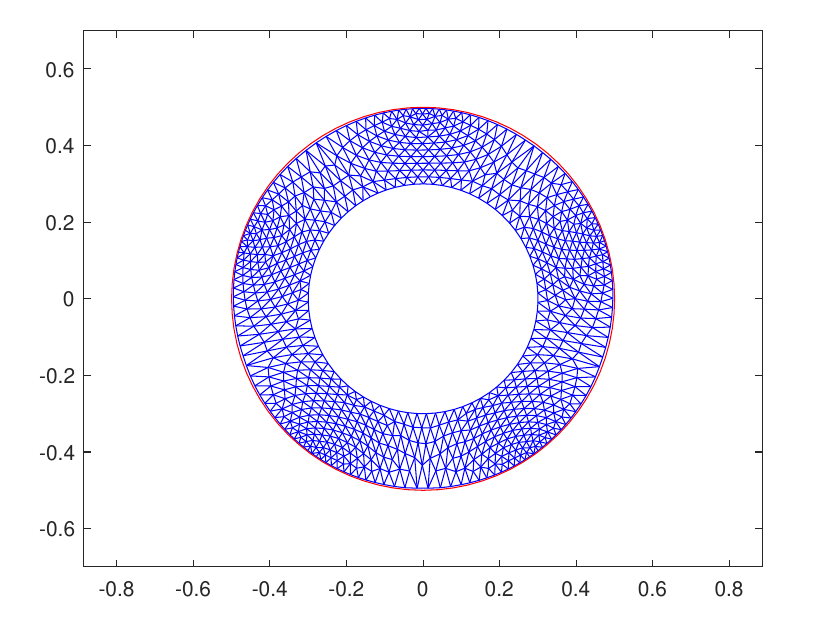}
}
\quad
\subfigure[$t=0.3$]{
\includegraphics[width=4.3cm]{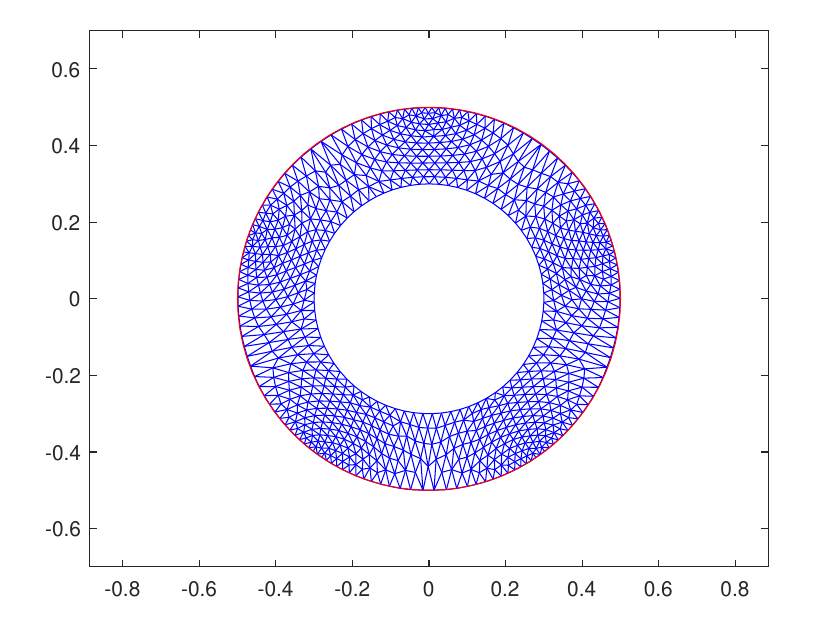}
}
\caption{Example~\ref{fbp-ex1}. The mesh of $N=1567$ is plotted at $t=0$, 0.01, 0.02, 0.05, 0.1, and 0.3 for $\lambda = -2/\ln(0.6)$.}
\label{fig:fbp-ex1-4}
\end{figure}

\begin{figure}[htbp]
	\centering
	\includegraphics[width=4.7cm]{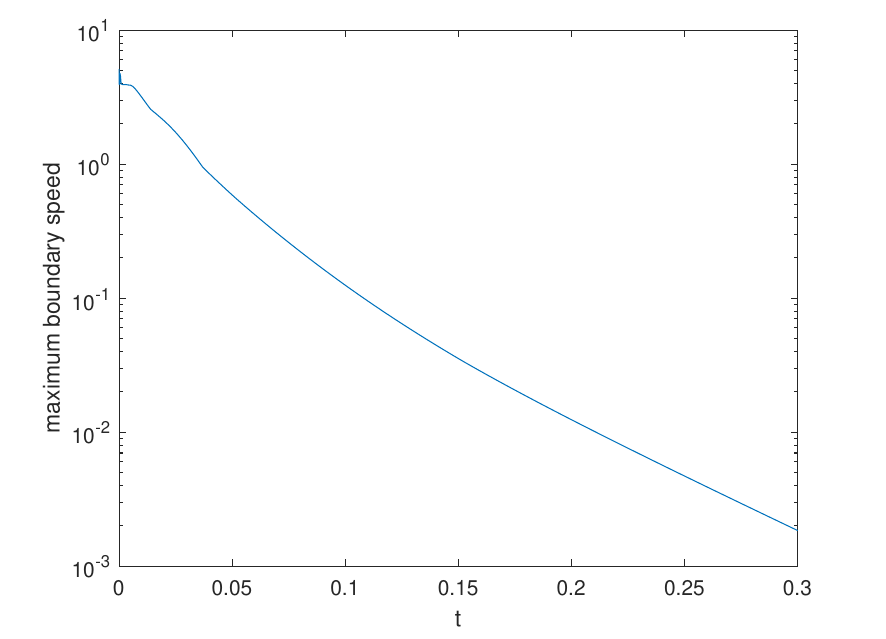}
	\caption{Example~\ref{fbp-ex1}. The maximum boundary velocity is plotted as a function of time
		for $\lambda = -2/\ln(0.6)$ and $N = 1567$.}
	\label{fig:fbp-ex1-5}
\end{figure}

\begin{exam}[\textbf{Exterior Bernoulli FBP - Accuracy test}]
\label{fbp-ex1}
This example is selected from Rabago \cite{Rabago2020}, where $\Gamma_1$ and the initial position of $\Gamma_2$ are taken as
the circles centered at the origin with radii
0.3 and 0.6, respectively, and $\lambda = - 2/ \ln (0.6)$. FBP (\ref{fbp-1}) has the exact solution $u = \ln (2 r)/\ln(0.6)$
and $\Gamma_2$ being the circle with radius 0.5. We compute the error as the average of the difference between
the radii of the boundary vertices on $\Gamma_2$ and the exact radius $0.5$ when the stopping criterion (toward steady state) is met.

A mesh at various time instants is plotted in Fig.~\ref{fig:fbp-ex1-2} and the corresponding maximum boundary velocity
is plotted as a function of time in Fig.~\ref{fig:fbp-ex1-3}. Notice that the maximum boundary velocity as a function of time
can be regarded as the convergence history towards the steady-state solution. From the figures we can see that
the maximum boundary velocity decreases gradually and the domain is converging towards steady state.
Fig.~\ref{fig:fbp-ex1-1} shows the convergence histories as the mesh is refined
for the error in the boundary location  for two strategies
of computing solution gradient used in boundary update.
The results show that the quadratic least squares fitting leads to second-order convergence whereas the area-weighted averaging
gives only first-order convergence.

We also consider a different initial position for $\Gamma_2$:
$x^2 + y^2 = (0.5+0.1 \sin(5 \arctan (y/x))$, to see how robust the moving mesh FEM is.
The mesh and maximum boundary velocity are shown in Figs.~\ref{fig:fbp-ex1-4} and
\ref{fig:fbp-ex1-5}, respectively. Once again, the results demonstrate the convergence towards steady state.
Interestingly, part of the initial position of $\Gamma_2$ is inside while the rest is outside
the exact solution circle (the circle with radius 0.5).  From Fig.~\ref{fig:fbp-ex1-4} we can see that the boundary vertices
initially inside the circle with radius 0.5 are moving outward and those outside the circle are moving inward, all towards
the exact solution circle. This is consistent with the formal analysis in Section~\ref{SEC:PTC} (also cf. Fig.~\ref{fig:boundary-movement}).
\qed
\end{exam}


\begin{figure}[htbp]
\centering
\subfigure[$t=0$]{
\includegraphics[width=4.3cm]{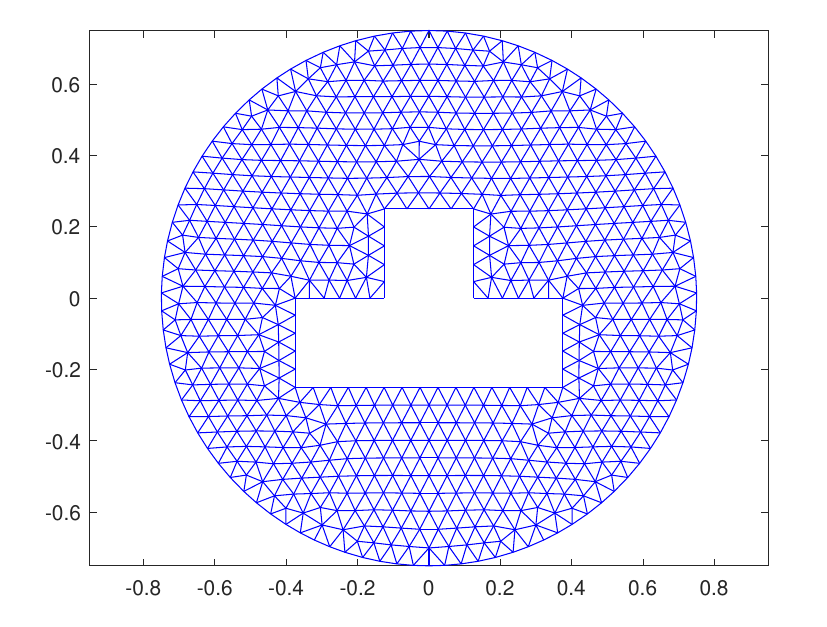}
}
\quad
\subfigure[$t=0.01$]{
\includegraphics[width=4.3cm]{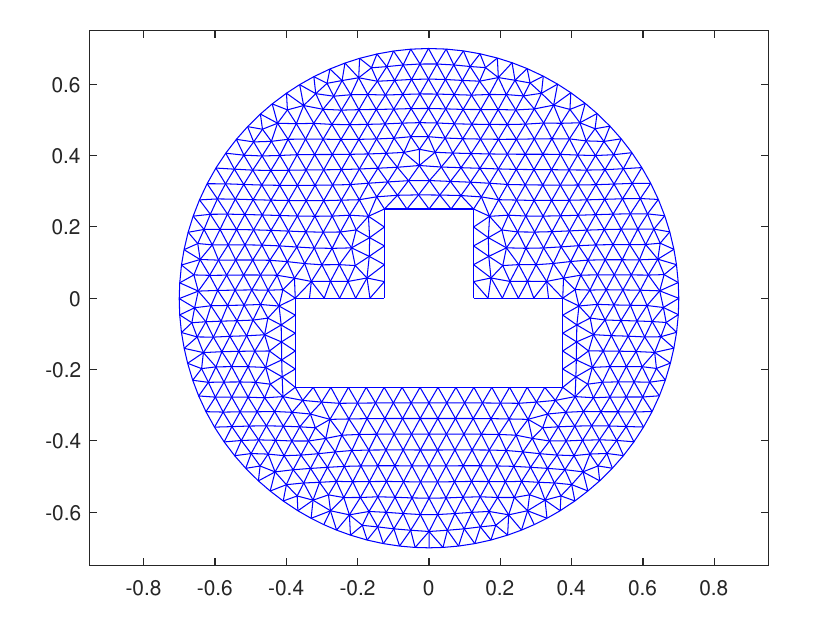}
}
\quad
\subfigure[$t=0.05$]{
\includegraphics[width=4.3cm]{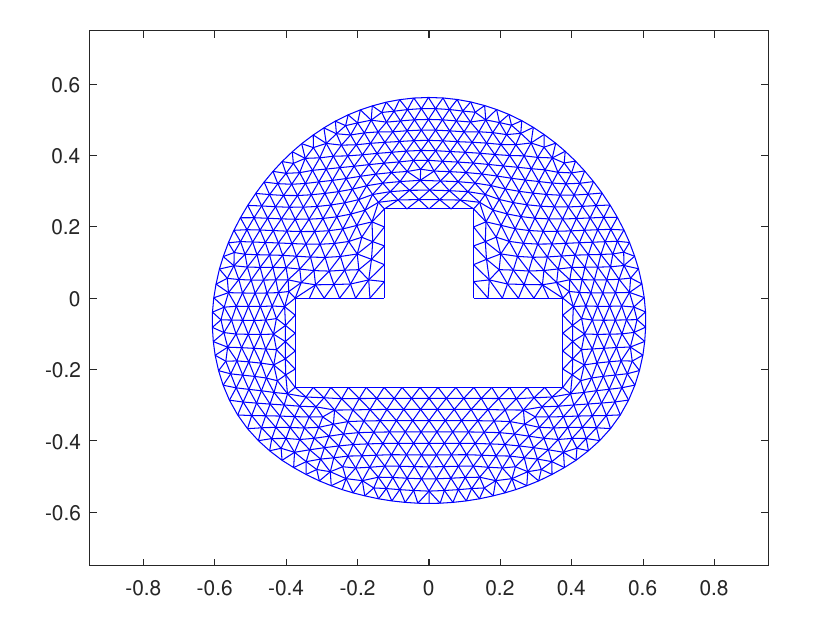}
}
\quad
\subfigure[$t=0.15$]{
\includegraphics[width=4.3cm]{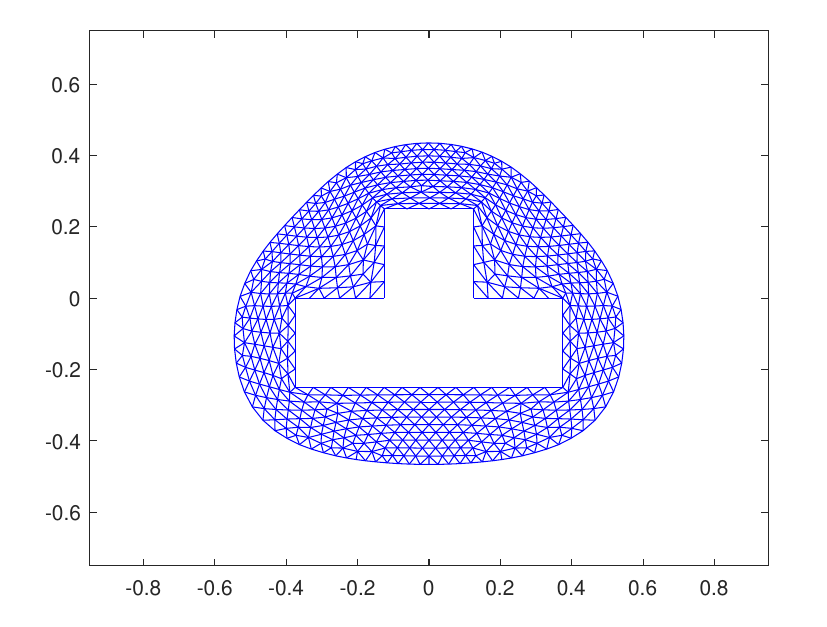}
}
\quad
\subfigure[$t=0.3$]{
\includegraphics[width=4.3cm]{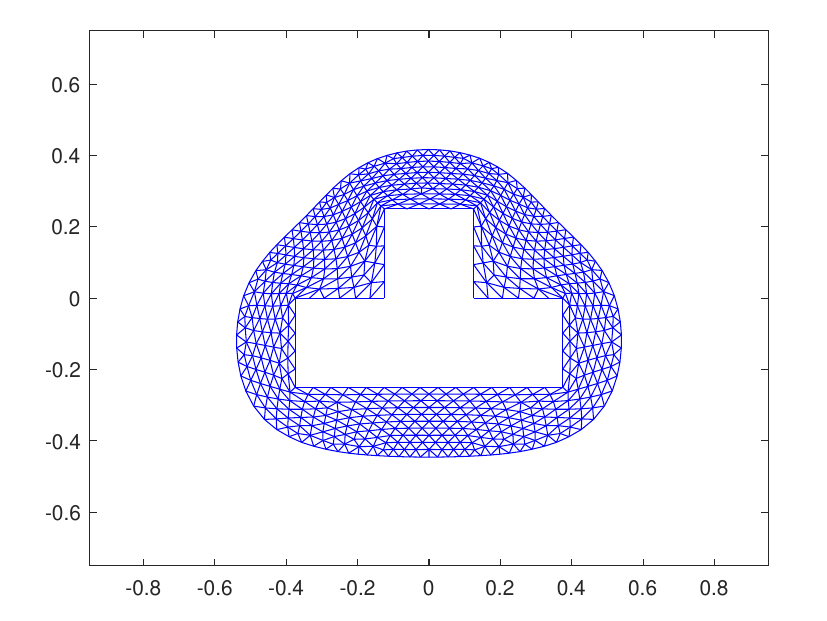}
}
\subfigure[$t=0.446$]{
\includegraphics[width=4.3cm]{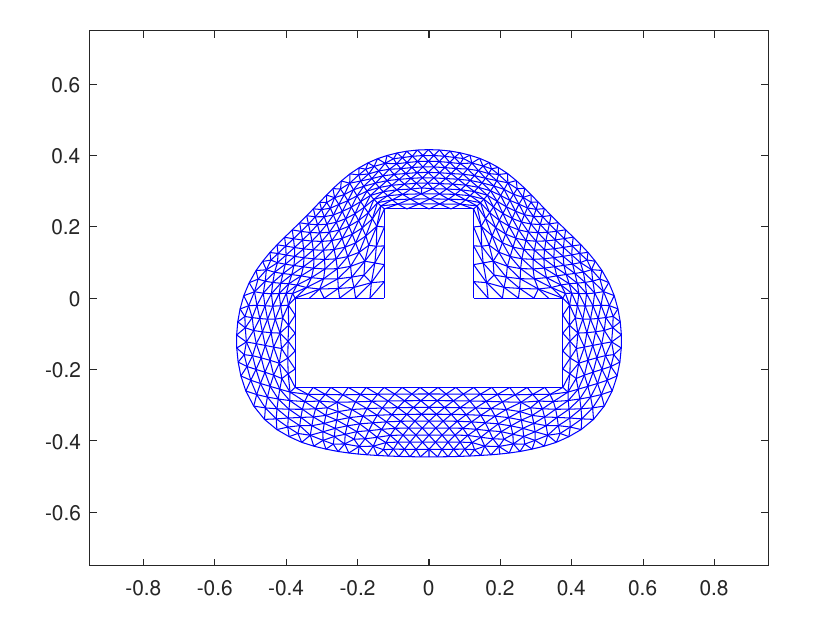}
}
\caption{Example~\ref{fbp-ex2}. The mesh of $N = 1259$ is plotted at $t=0$, 0.05, 0.1, 0.15, 0.3, and 0.446 for $\lambda = 5$.}
\label{fig:fbp-ex2-1}
\end{figure}

\begin{figure}[htbp]
\centering
\includegraphics[width=4.7cm]{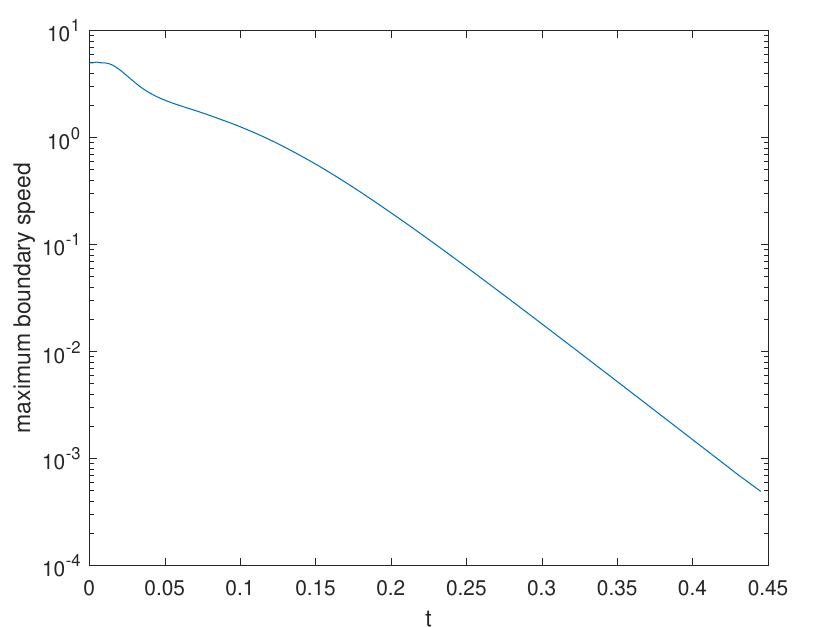}
\caption{Example~\ref{fbp-ex2}. The maximum boundary velocity is plotted as a function of time
for $\lambda = 5$ and $N = 1259$.}
\label{fig:fbp-ex2-2}
\end{figure}

\begin{figure}[htbp]
\centering
\includegraphics[width=4.7cm]{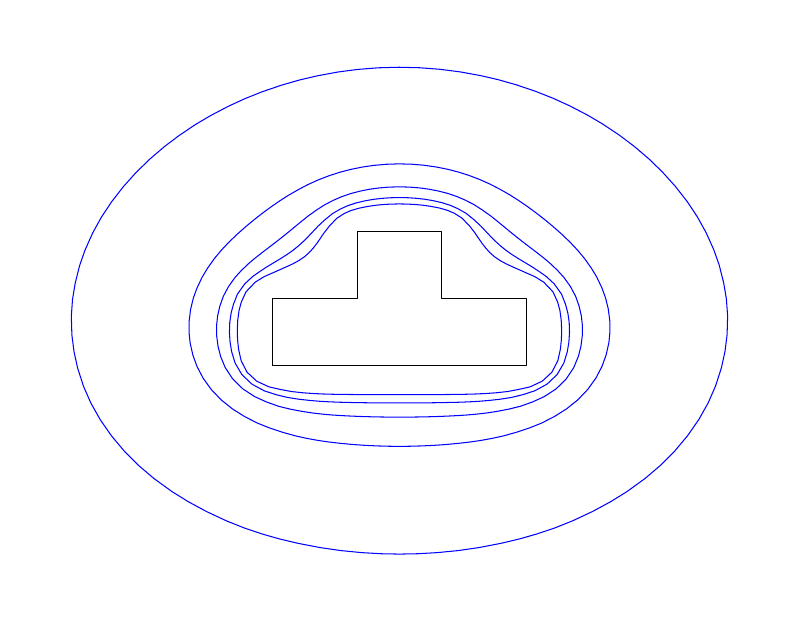}
\caption{Example~\ref{fbp-ex2}. $\Gamma_2$ obtained for $\lambda=1$, 3, 5, 7, and 9 with a mesh of $N = 1259$.}
\label{fig:fbp-ex2-3}
\end{figure}

\begin{exam}[\textbf{Exterior Bernoulli FBP with $T$-shape}]
\label{fbp-ex2}
For this example, $\Gamma_1$ is taken as the boundary of the $T$-shape
\[
(-3/8, 3/8) \times (-1/4, 0) \cup (-1/8, 1/8) \times [0, 1/4)
\]
and the initial position of $\Gamma_2$ is a circle of radius 0.75. This problem was used by several researchers
(e.g., see Eppler and Harbrecht \cite{Eppler2012}).

Fig.~\ref{fig:fbp-ex2-1} shows the mesh at various time instants for $\lambda = 5$. Fig.~\ref{fig:fbp-ex2-2} shows that
the maximum boundary velocity decreases as  the time increases, implying that (\ref{fbp-0}) has a steady-state solution for this example.
Fig.~\ref{fig:fbp-ex2-3} shows $\Gamma_2$ obtained for $\lambda=1$, 3, 5, 7 and 9. As $\lambda$ increases,
$\Gamma_2$ is getting closer to $\Gamma_1$. The results obtained here are comparable
with those in literature and particularly those obtained in \cite{Eppler2012} using a shape optimization method.
\qed
\end{exam}

\begin{figure}[htbp]
\centering
\subfigure[$t=0$]{
\includegraphics[width=4.3cm]{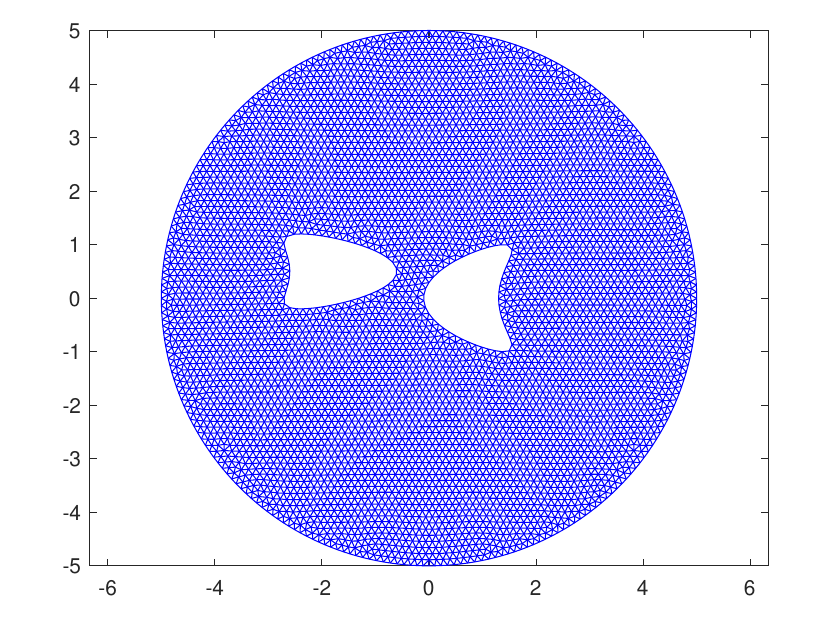}
}
\quad
\subfigure[$t=1$]{
\includegraphics[width=4.3cm]{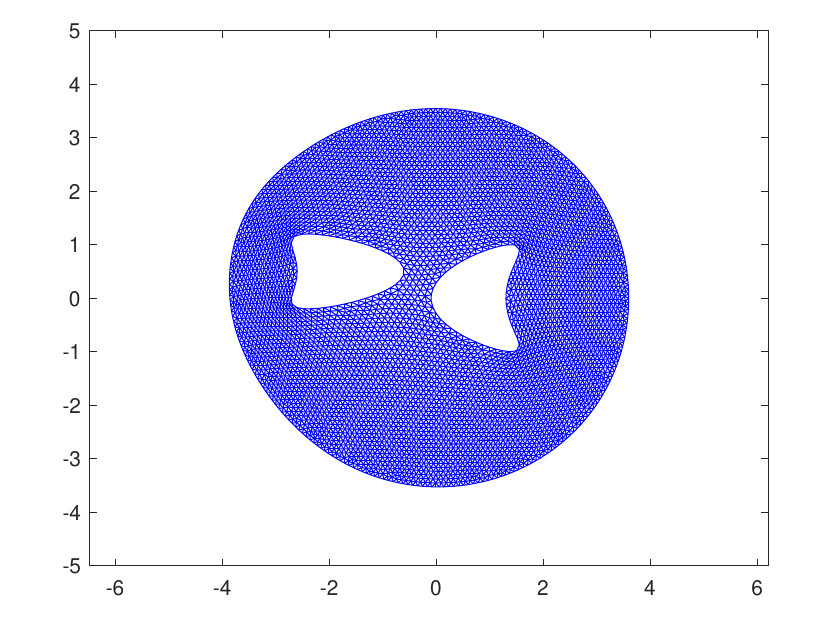}
}
\quad
\subfigure[$t=2$]{
\includegraphics[width=4.3cm]{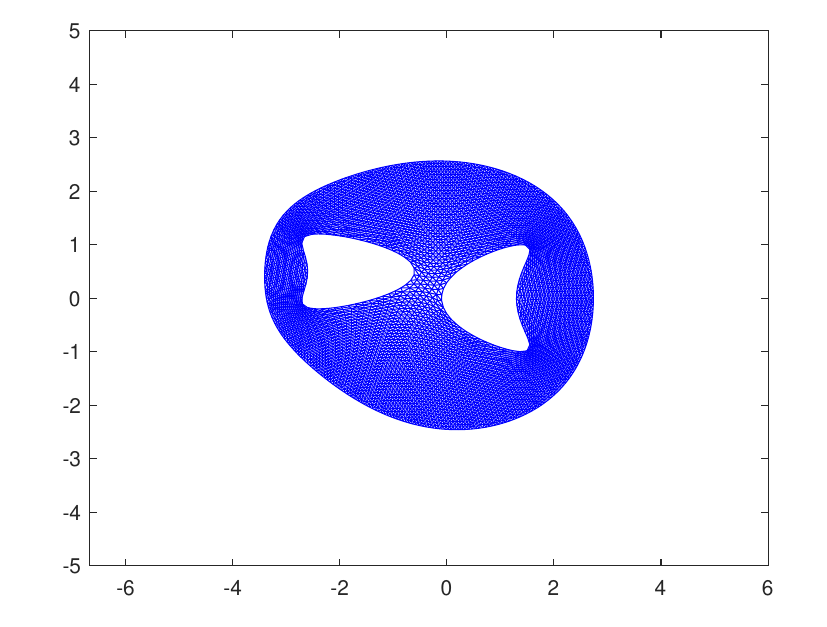}
}
\quad
\subfigure[$t=3$]{
\includegraphics[width=4.3cm]{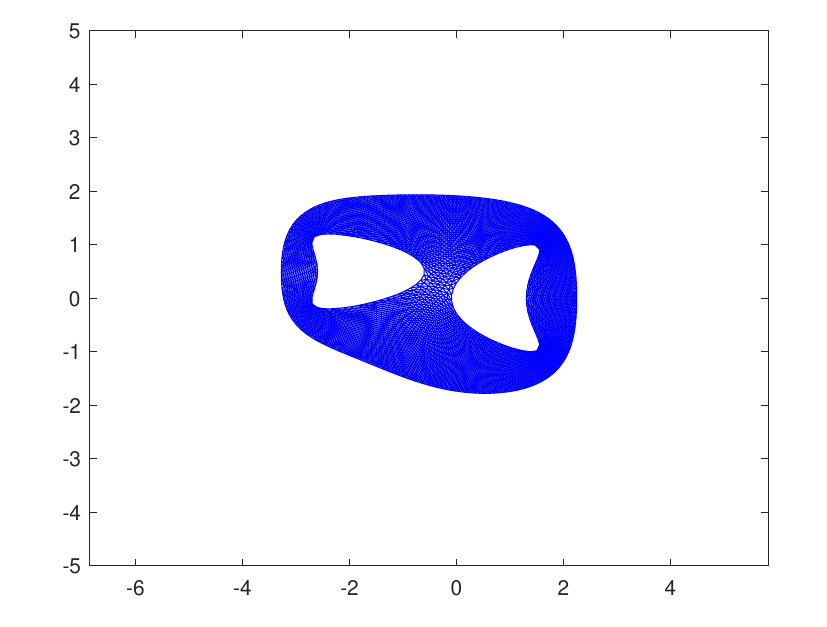}
}
\quad
\subfigure[$t=6$]{
\includegraphics[width=4.3cm]{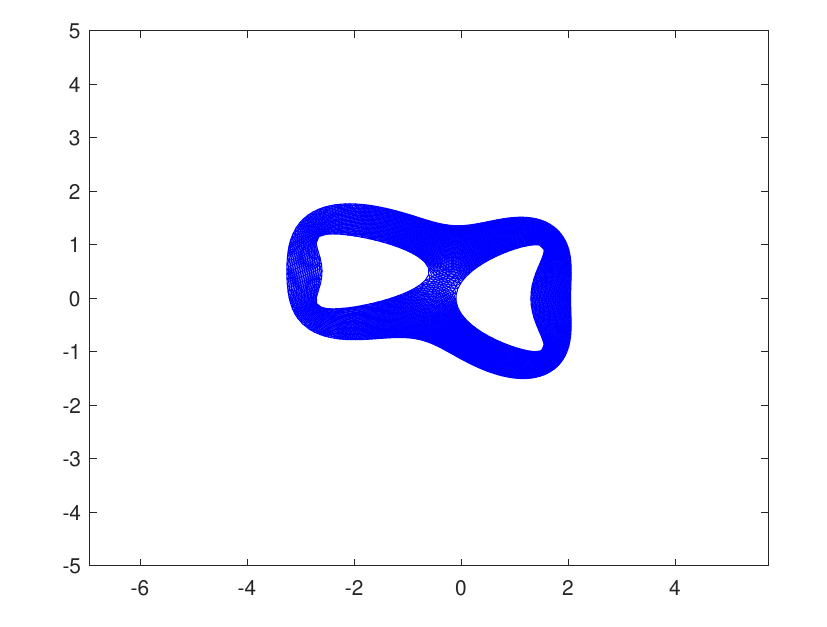}
}
\quad
\subfigure[$t=9.941$]{
\includegraphics[width=4.3cm]{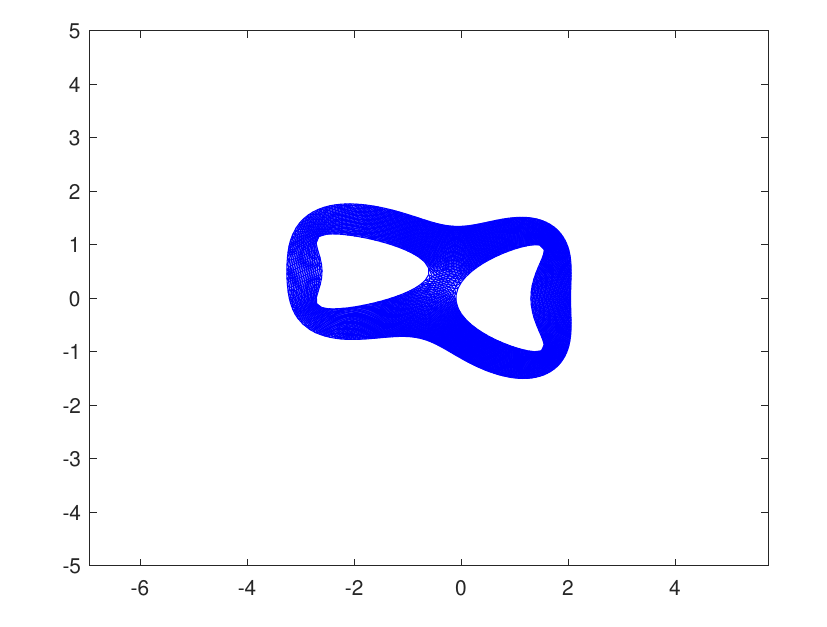}
}
\caption{Example~\ref{fbp-ex3}. The mesh of $N=10630$ is plotted at $t=0$, 1, 2, 3, 6, and 9.941
for $\lambda = 1.5$.}
\label{fig:fbp-ex3-1}
\end{figure}

\begin{figure}[ht!]
\centering
\includegraphics[width=4.7cm]{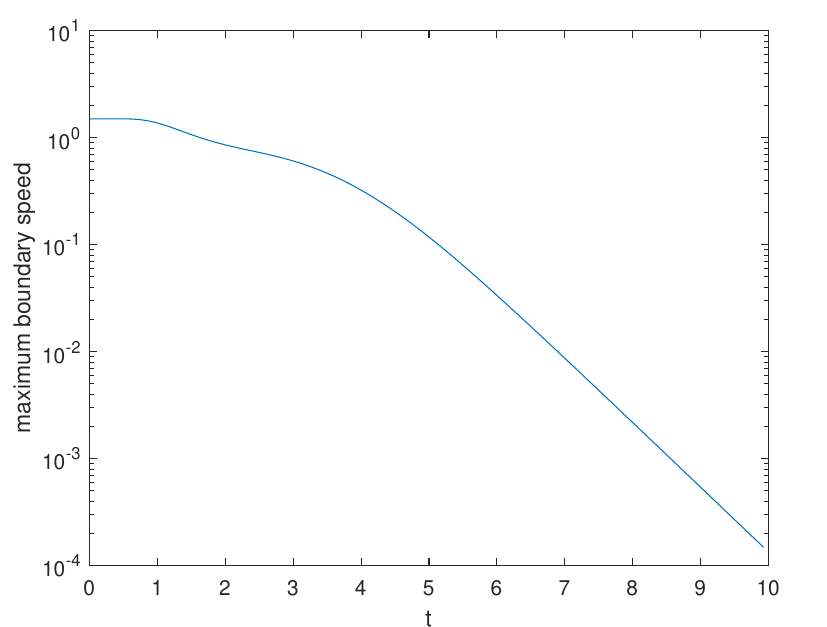}
\caption{Example~\ref{fbp-ex3}. The maximum boundary velocity is plotted as a function of time
for $\lambda = 1.5$ and $N = 10630$.}
\label{fig:fbp-ex3-2}
\end{figure}

\begin{figure}[ht!]
	\centering
	\includegraphics[width=5cm]{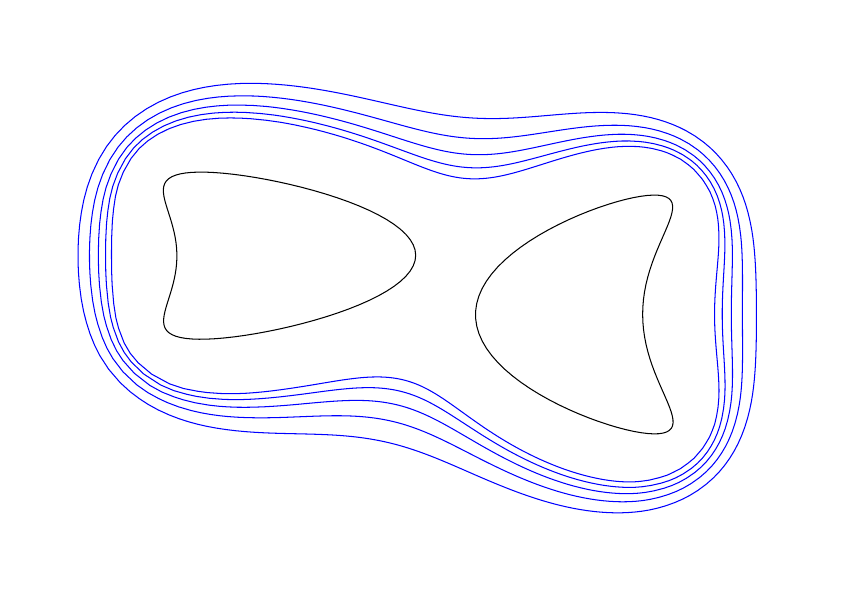}
	\caption{Example~\ref{fbp-ex4}. The boundary $\Gamma_2$ is obtained with $\lambda= 1.1$, 1.3, 1.5, 1.7 and 1.9
		and a mesh of $N = 10630$.}
	\label{fig:fbp-ex3-3}
\end{figure}

\begin{exam}[\textbf{Exterior Bernoulli FBP with two disjoint shapes}]
\label{fbp-ex3}

This example is selected from Rabago \cite{Rabago2020}. The interior boundary, $\Gamma_1$,
consists of the boundary of two disjoint shapes
\begin{align*}
& (1+0.7\cos(\theta)-0.4 \cos(2\theta), \; \sin(\theta)),\; 0\le \theta \le 2 \pi
\\
& (-2+\cos(\theta)+0.4\cos(2 \theta), \; 0.5+0.7\sin(\theta)), \; 0\le \theta \le 2 \pi .
\end{align*}
The initial position of $\Gamma_2$ is taken as a circle of radius 5 with center (0,0).

Figs.~\ref{fig:fbp-ex3-1} and \ref{fig:fbp-ex3-2} show a mesh at various time instants and
the maximum boundary velocity as a function of time, respectively.
The location of $\Gamma_2$ obtained for several values of $\lambda$ is plotted in Fig.~\ref{fig:fbp-ex3-3}.
The results show that the moving mesh FEM with the pseudo-transient continuation
works well for this example with more complex $\Gamma_1$.
Particularly, the mesh stays free of tangling for all computations.
\qed
\end{exam}


\begin{figure}[htbp]
\centering
\subfigure[$t=0$]{
\includegraphics[width=4.3cm]{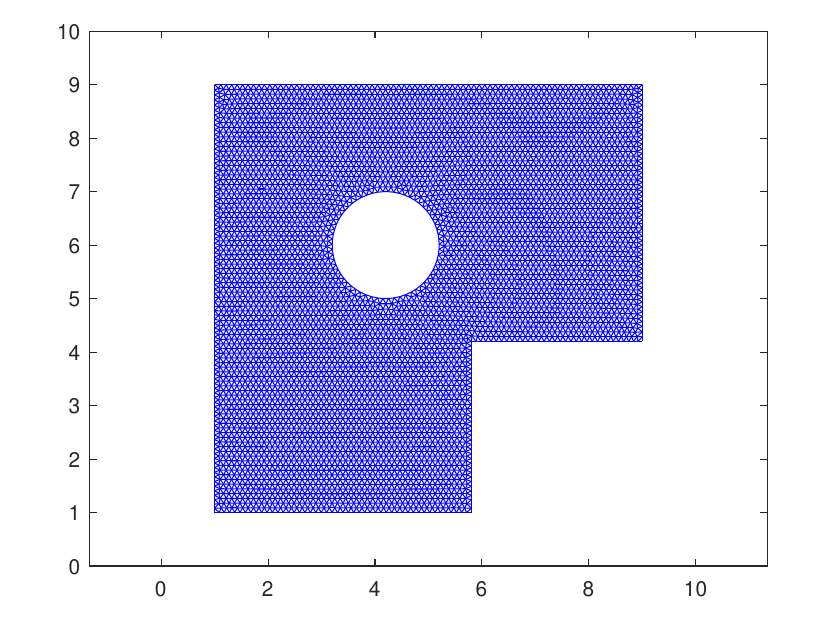}
}
\quad
\subfigure[$t=4$]{
\includegraphics[width=4.3cm]{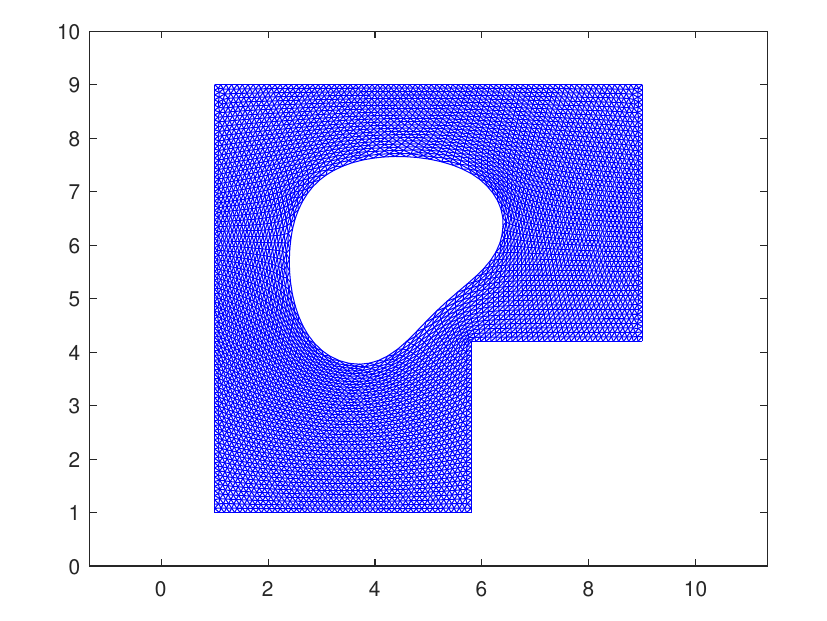}
}
\quad
\subfigure[$t=8$]{
\includegraphics[width=4.3cm]{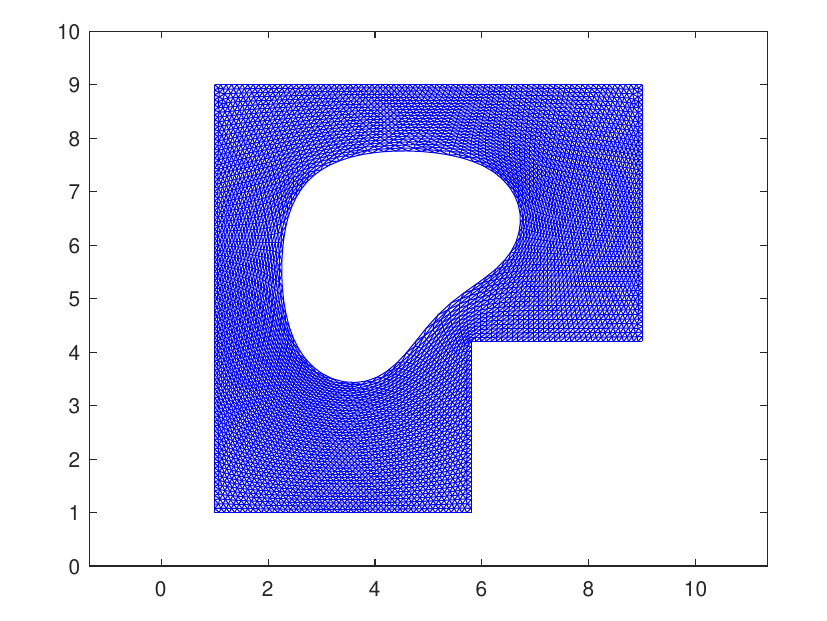}
}
\quad
\subfigure[$t=10$]{
\includegraphics[width=4.3cm]{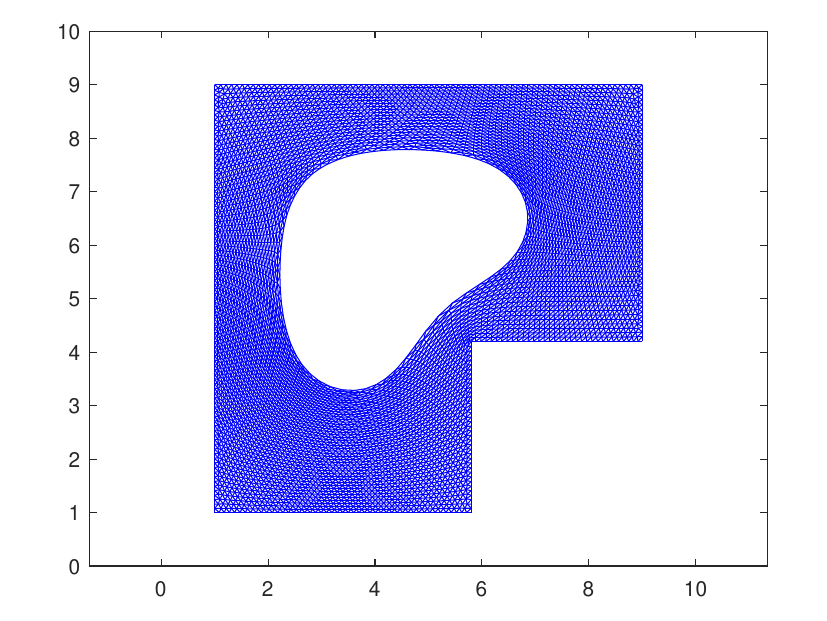}
}
\quad
\subfigure[$t=16$]{
\includegraphics[width=4.3cm]{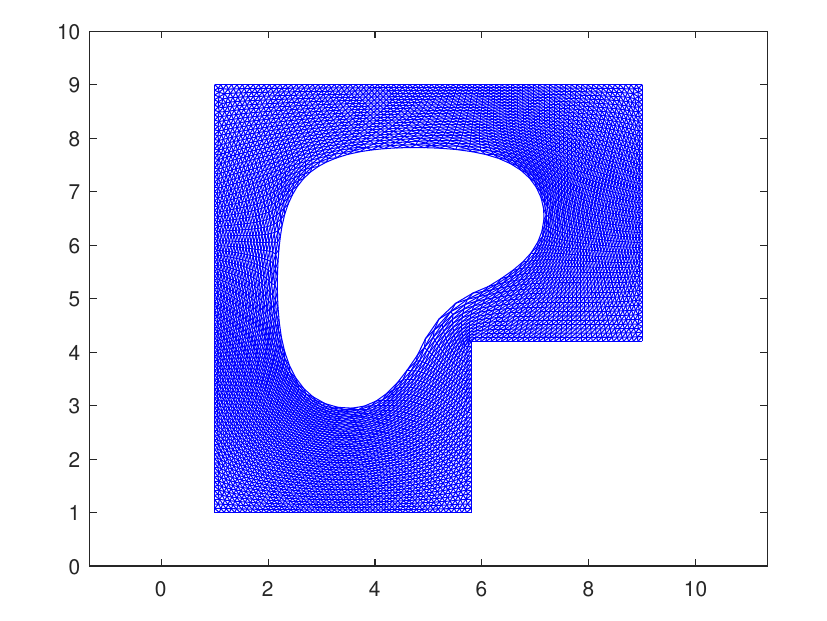}
}
\quad
\subfigure[$t=20$]{
\includegraphics[width=4.3cm]{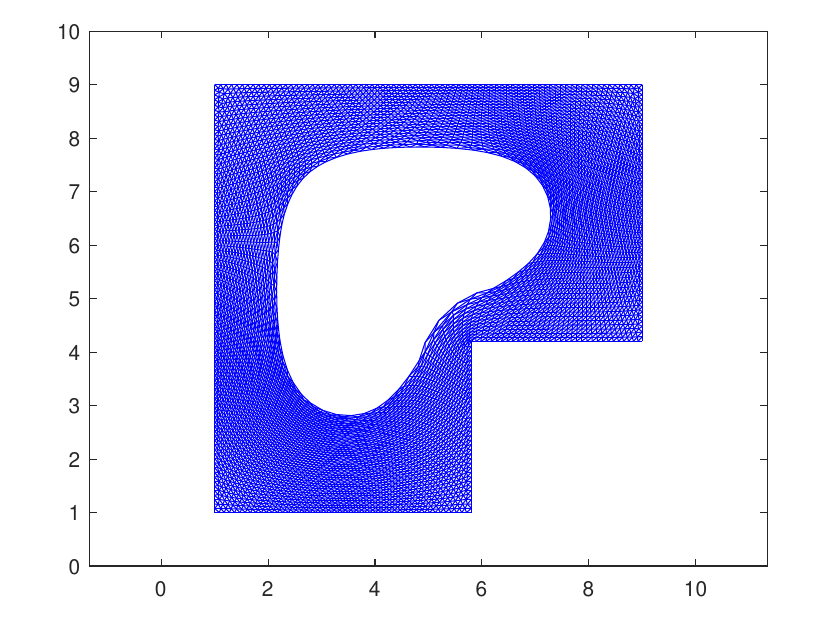}
}
\caption{Example~\ref{fbp-ex4}. The mesh of $N=11334$ is plotted at $t=0$, 4, 8, 10, 16, and 20 for $\lambda = 0.9$.}
\label{fig:fbp-ex4-1}
\end{figure}

\begin{figure}[htbp!]
\centering
\includegraphics[width=4.7cm]{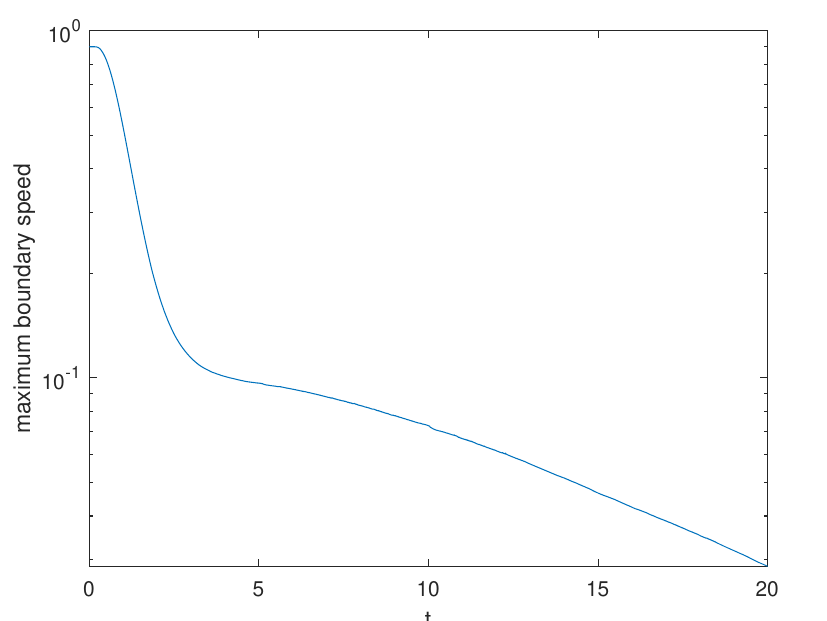}
\caption{Example~\ref{fbp-ex4}. The maximum boundary velocity is plotted as a function of time for $\lambda = 0.9$.}
\label{fig:fbp-ex4-2}
\end{figure}

\begin{figure}[ht!]
\centering
\includegraphics[width=5cm]{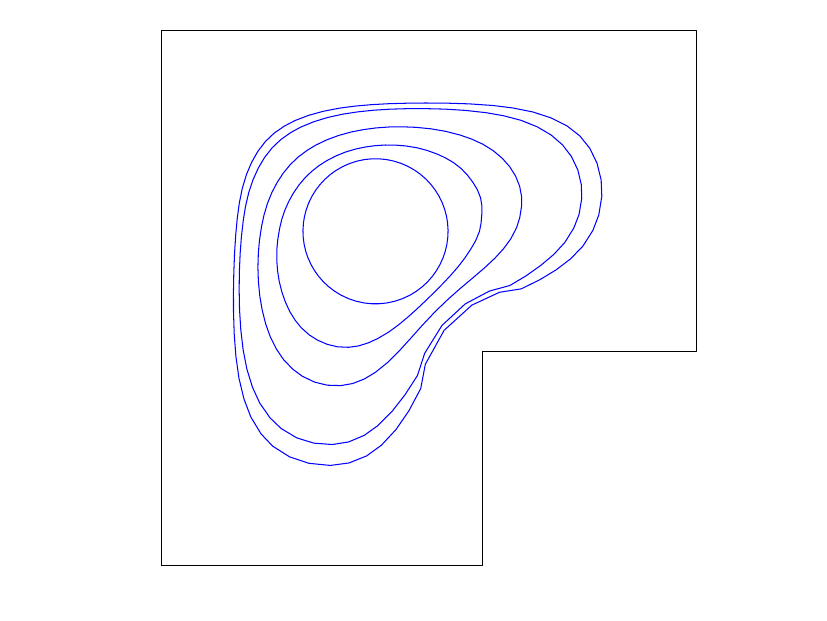}
\caption{Example~\ref{fbp-ex4}. The boundary $\Gamma_2$ is obtained with $\lambda= 0.75$, 0.8, 0.85, 0.9 and 0.95.}
\label{fig:fbp-ex4-3}
\end{figure}

\begin{exam}[\textbf{Interior Bernoulli FBP with $L$-shape}]
\label{fbp-ex4}

Finally, we consider an interior Bernoulli FBP. In this example, $\Gamma_1$ is taken as the boundary of the $L$-shape
\[
(1, 5.8) \times (1, 9) \cup [5.8, 9) \times (4.2, 9)
\]
and the initial position of $\Gamma_2$ is the circle of radius 1.5 with center (4.2, 6). A similar example was considered
by Flucher and Rumpf \cite{Flucher1997} and several other researchers.

Fig.~\ref{fig:fbp-ex4-1} shows the mesh at various time instants and Fig.~\ref{fig:fbp-ex4-2} shows the maximum boundary velocity
as a function of time. For this example, the solution converges more slowly to steady state than previous examples.
The computation is stopped at $t = 20$ when the maximum boundary velocity is about $2 \times 10^{-2}$ and
the boundary displacement is about $2\times 10^{-5}$. Nevertheless, the figures show that the maximum boundary velocity
decreases steadily and $\Gamma_2$ is converging towards steady state.
Fig.~\ref{fig:fbp-ex4-3} shows $\Gamma_2$ for several values of $\lambda$. As $\lambda$ increases, $\Gamma_2$
is getting closer to $\Gamma_1$.
The results are comparable with those in \cite{Flucher1997} where the explicit and implicit Neumann methods are used.
\qed
\end{exam}


\pagebreak
\section{Numerical examples for FBPs with non-constant Bernoulli condition and nonlinear FBPs}
\label{SEC:numerics-2}

The moving mesh method described in Section~\ref{SEC:MM-FEM} can be used for more general FBPs without major modifications.
To demonstrate this, we present in this section numerical results for three examples, one with non-constant Bernoulli
boundary condition, one with  the $p$-Laplacian (nonlinear), and one being a nonlinear obstacle problem.
FBPs with non-constant Bernoulli conditions and/or $p$-Laplacian have been studied by a number of researchers,
e.g., see Acker and Meyer \cite{Acker-1995} and Henrot and Shahgholian \cite{Henrot-2002}.
Obstacle problems are a classical and important types of FBPs (e.g., see Ros-Oton \cite{Ros-Oton-2018}).
The settings and values of the parameters used in the computation are the same as in the previous section.

\begin{figure}[htbp]
\centering
\subfigure[$t=0$]{
\includegraphics[width=4.3cm]{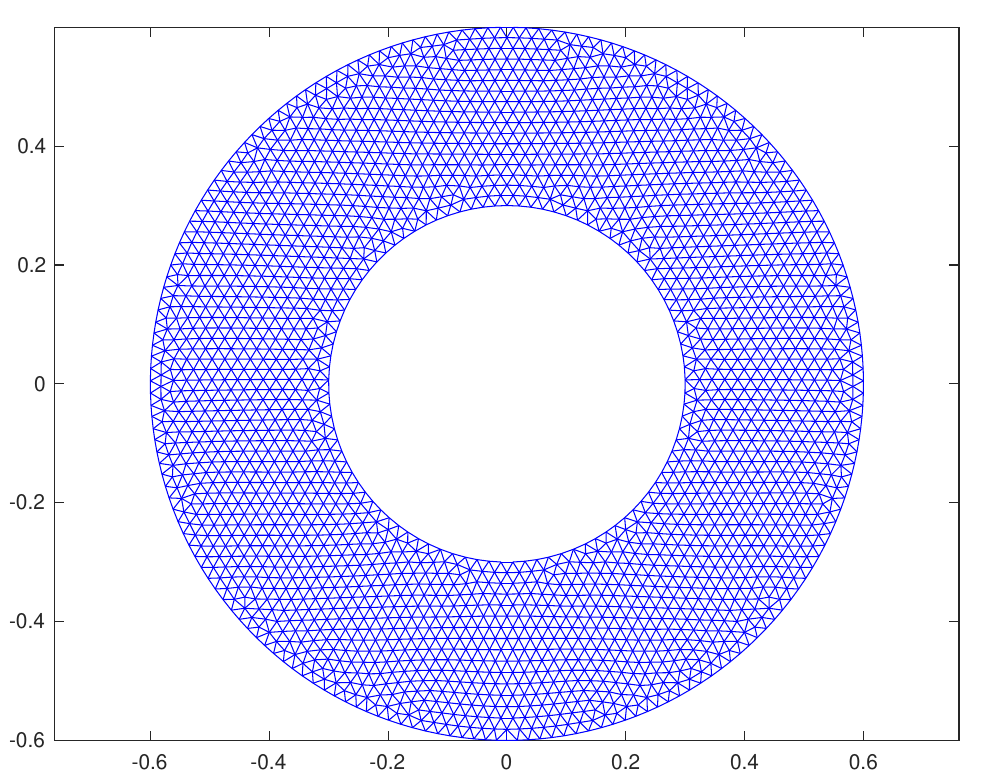}
}
\quad
\subfigure[$t=0.15$]{
\includegraphics[width=4.3cm]{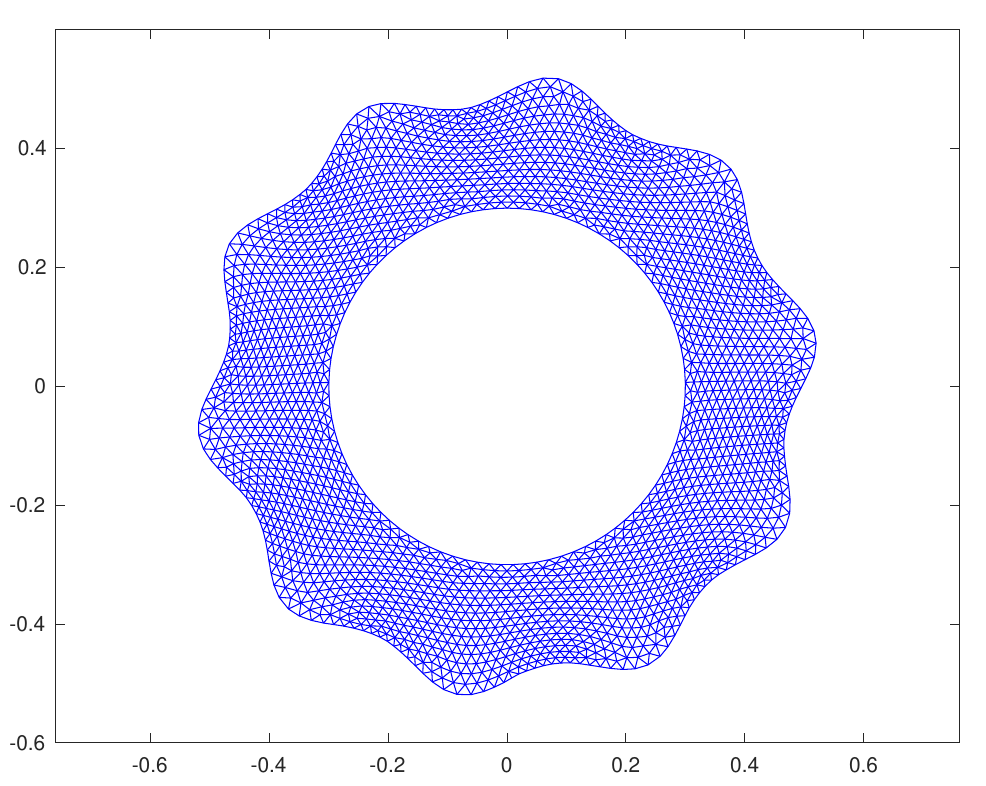}
}
\\
\quad
\subfigure[$t=0.3$]{
\includegraphics[width=4.3cm]{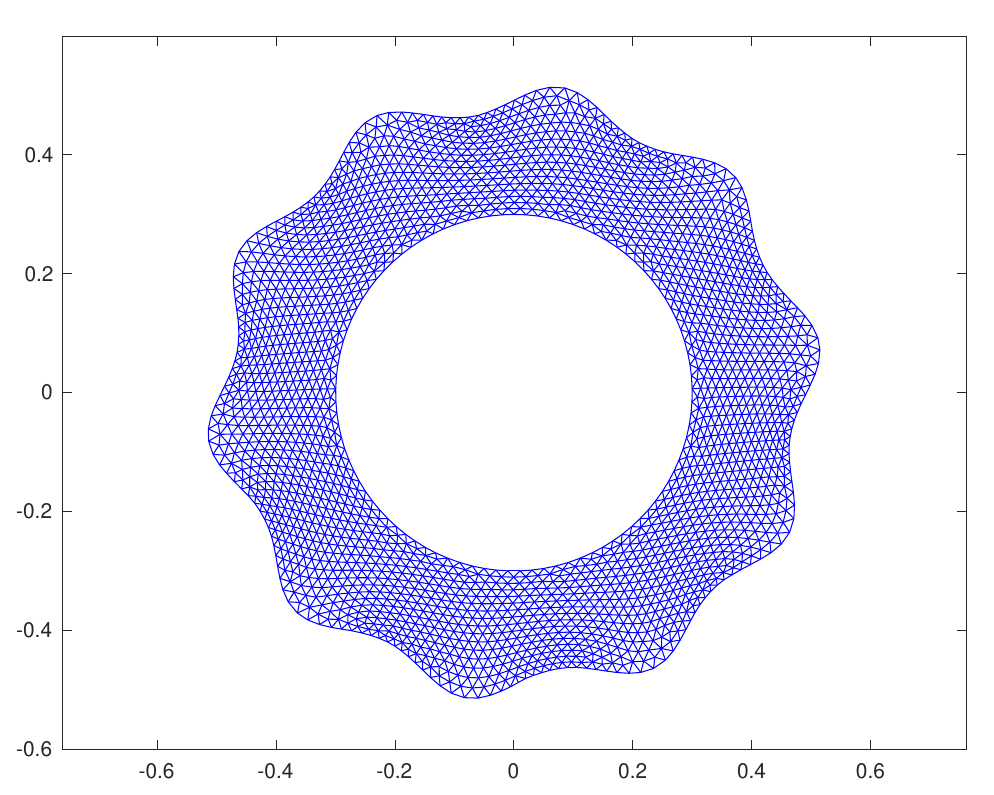}
}
\quad
\subfigure[$t=0.468$]{
\includegraphics[width=4.3cm]{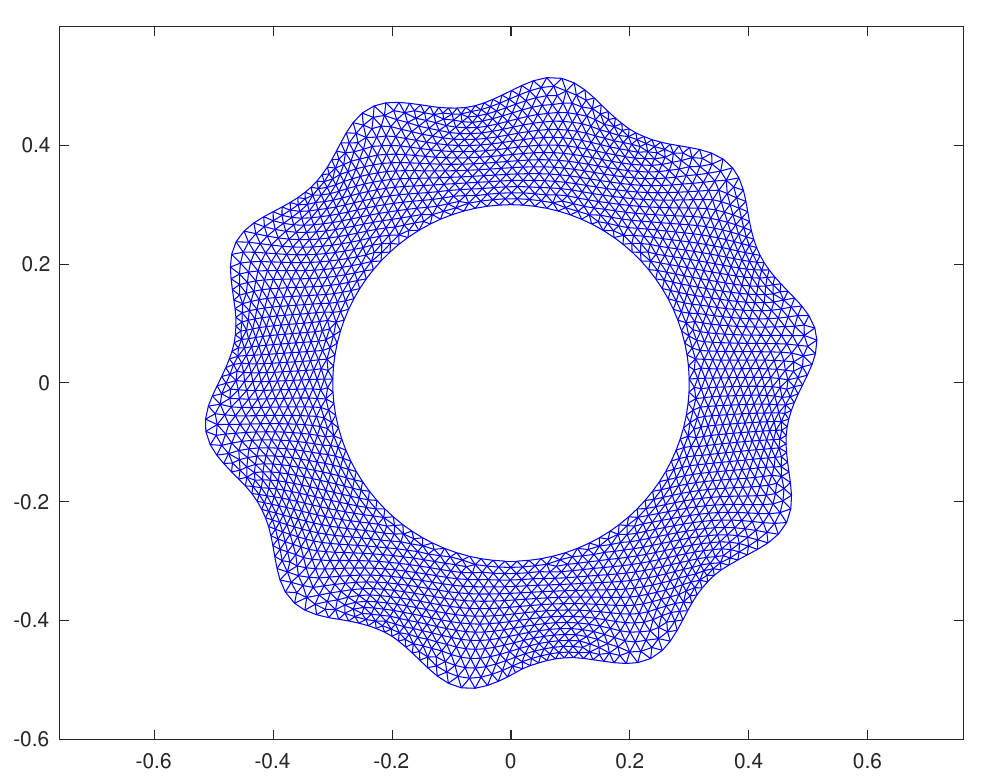}
}
\caption{Example~\ref{fbp-ex5}. The mesh of $N = 4618$ is plotted at $t=0$, 0.15, 0.3, and 0.468
for variable $\lambda$ (\ref{lambda-var-1}).}
\label{fig:fbp-ex5-1}
\end{figure}

\begin{figure}[htbp]
\centering
\includegraphics[width=4.7cm]{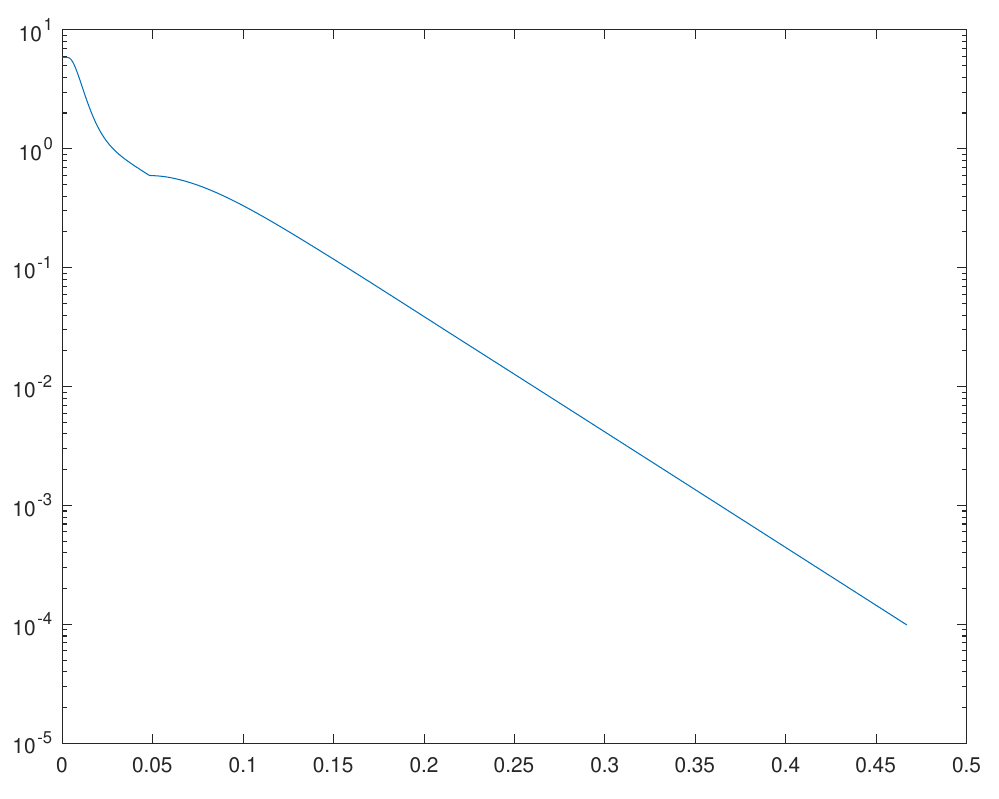}
\caption{Example~\ref{fbp-ex5}. The maximum boundary velocity is plotted as a function of time
for  variable $\lambda$ (\ref{lambda-var-1}) and $N =4618$.}
\label{fig:fbp-ex5-2}
\end{figure}

\begin{exam}[\textbf{Exterior Bernoulli FBP with non-constant Bernoulli condition}]
\label{fbp-ex5}
This example is the same as Example~\ref{fbp-ex1} except that a non-constant Bernoulli boundary condition is used,
\begin{equation}
\lambda = - \frac{2}{\ln (0.6)} \left ( 1 - 0.5\sin(10 \arctan(\frac{y}{x}) \right ) .
\label{lambda-var-1}
\end{equation}
The initial position of $\Gamma_2$ is taken as the circle with radius 0.6.
A mesh and the corresponding maximum boundary velocity are plotted in Figs.~\ref{fig:fbp-ex5-1}
and \ref{fig:fbp-ex5-2}, respectively. One can see that the steady-state $\Gamma_2$ for this example is a wavy circle,
which is different from a circle in Example~\ref{fbp-ex1}.
The results also show that the moving mesh FEM with the pseudo-transient continuation works well for this example.
\qed
\end{exam}


\begin{figure}[ht!]
\centering
\subfigure[$t=0$]{
\includegraphics[width=4.3cm]{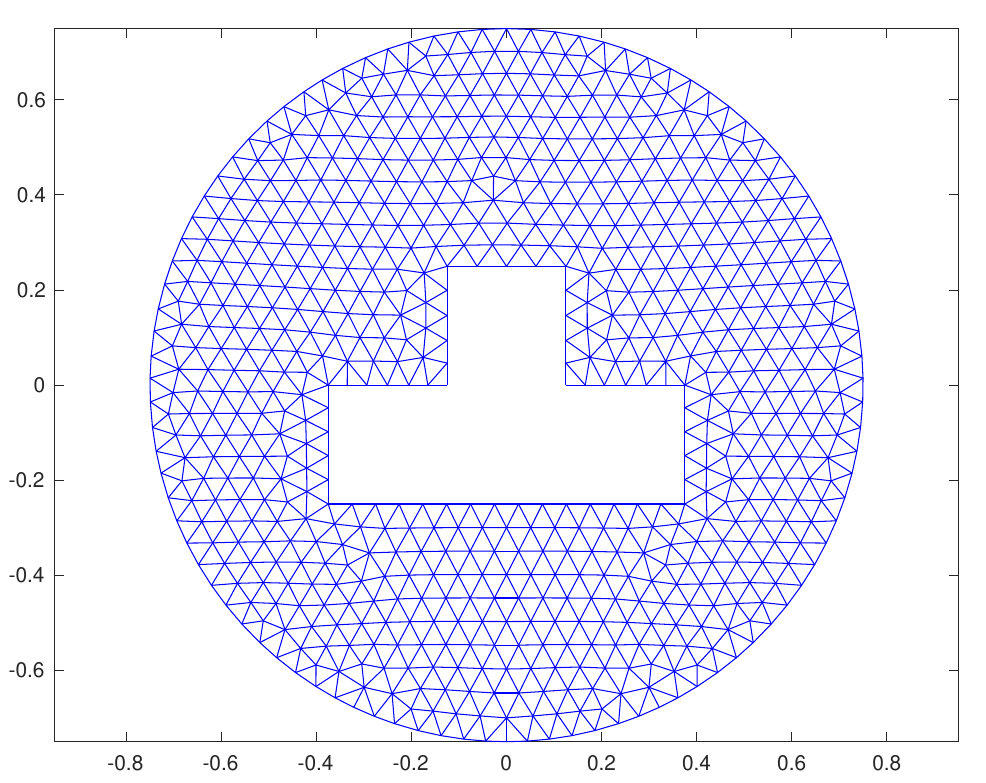}
}
\quad
\subfigure[$t=0.01$]{
\includegraphics[width=4.3cm]{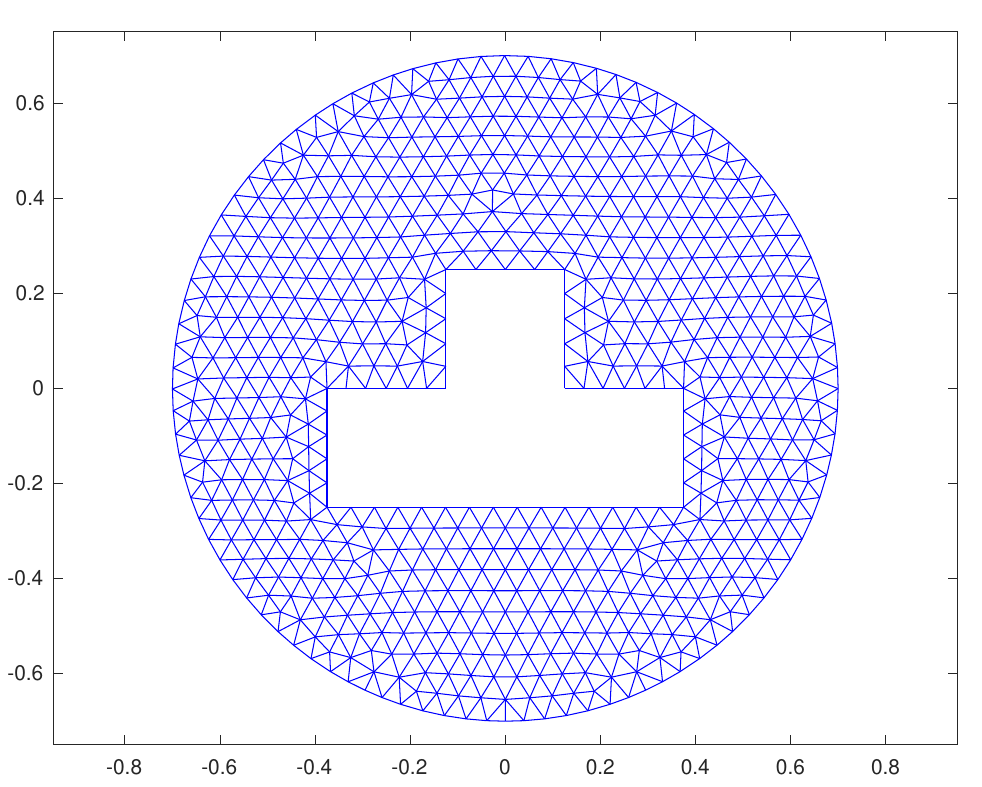}
}
\quad
\subfigure[$t=0.05$]{
\includegraphics[width=4.3cm]{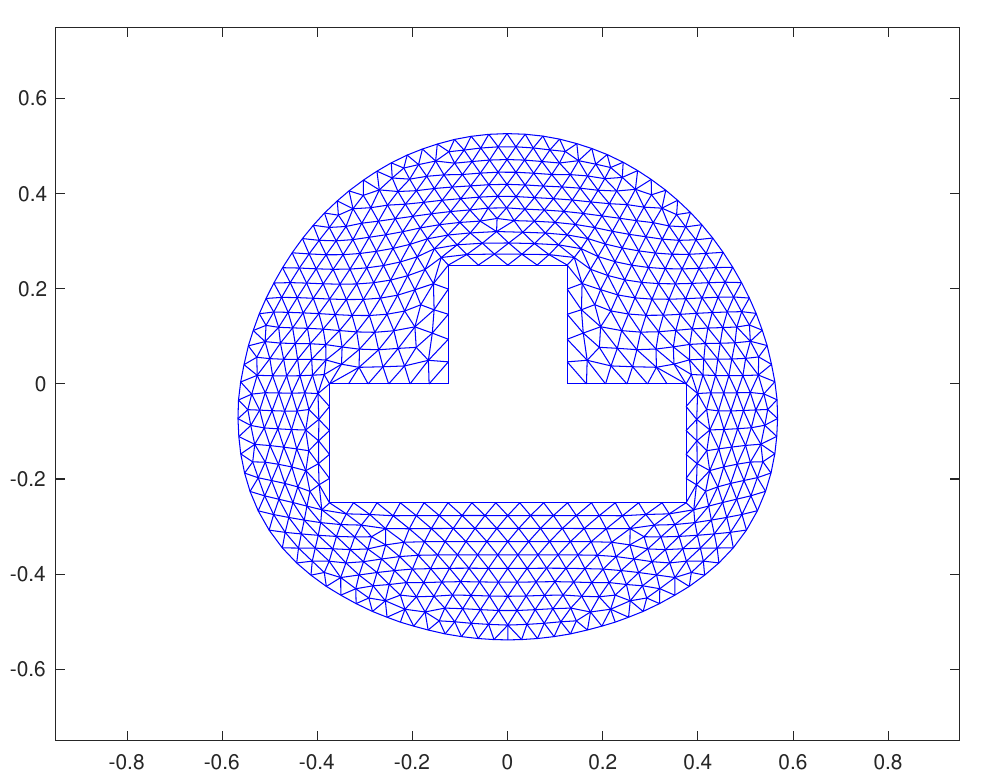}
}
\quad
\subfigure[$t=0.15$]{
\includegraphics[width=4.3cm]{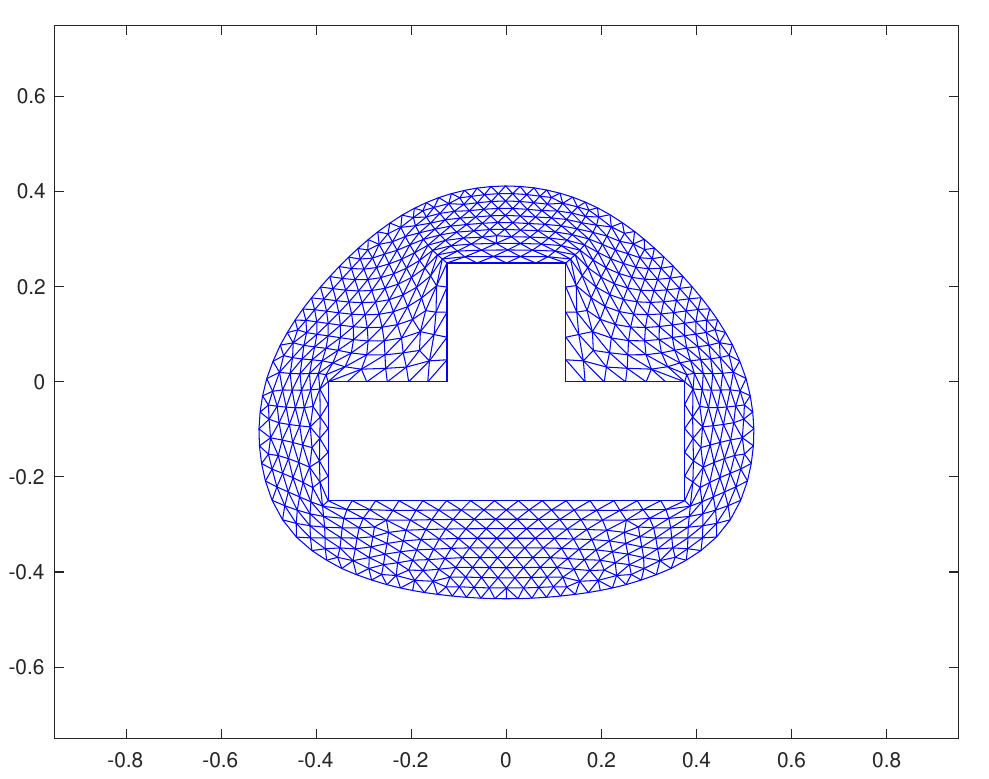}
}
\quad
\subfigure[$t=0.3$]{
\includegraphics[width=4.3cm]{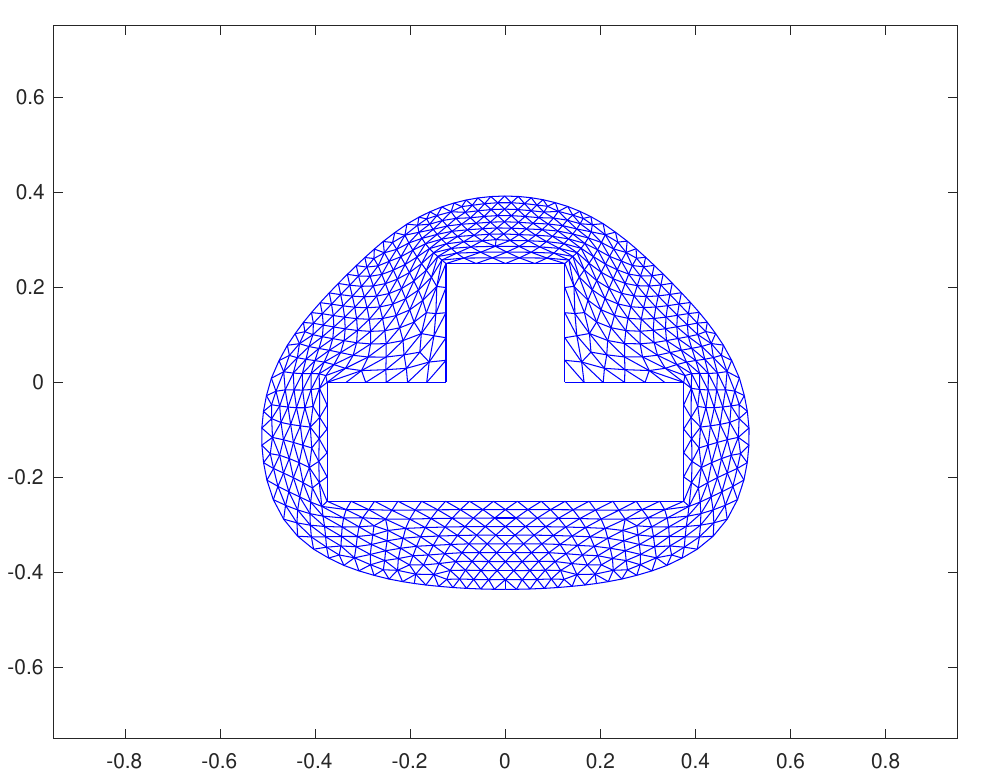}
}
\quad
\subfigure[$t=0.677$]{
\includegraphics[width=4.3cm]{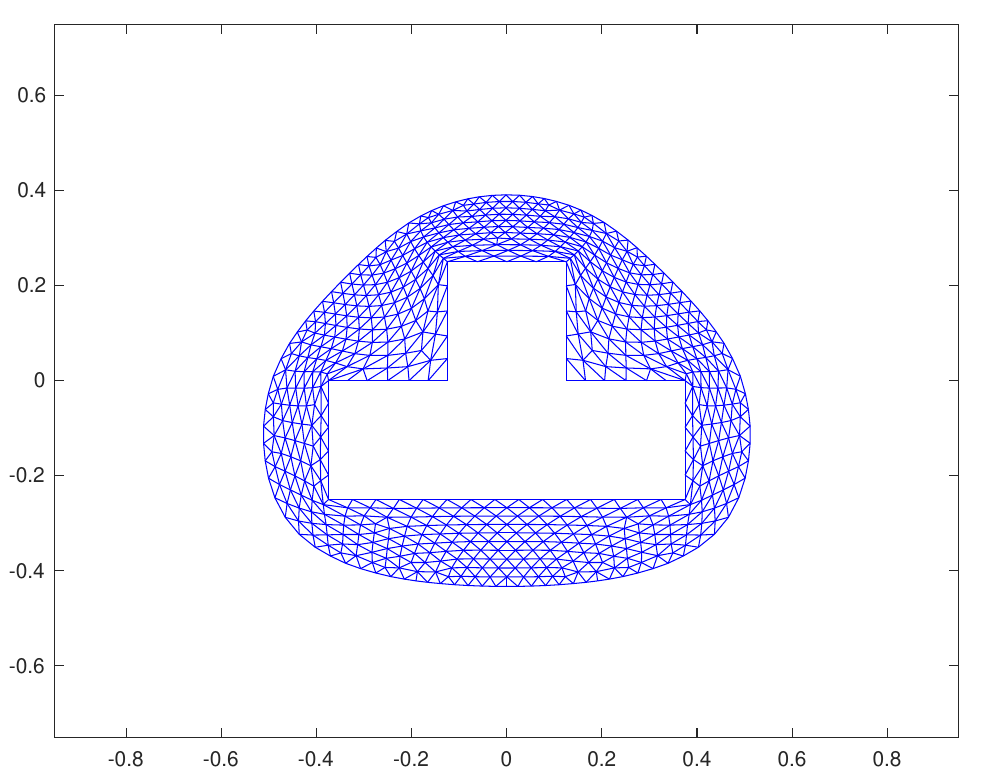}
}
\caption{Example~\ref{fbp-ex6} with $p = 1.5$. The mesh of $N = 1259$ is plotted at $t=0$, 0.05, 0.1, 0.15, 0.3, and 0.677 for $\lambda = 5$.}
\label{fig:fbp-ex6-1}
\end{figure}

\begin{figure}[htbp!]
\centering
\subfigure[$t=0$]{
\includegraphics[width=4.5cm]{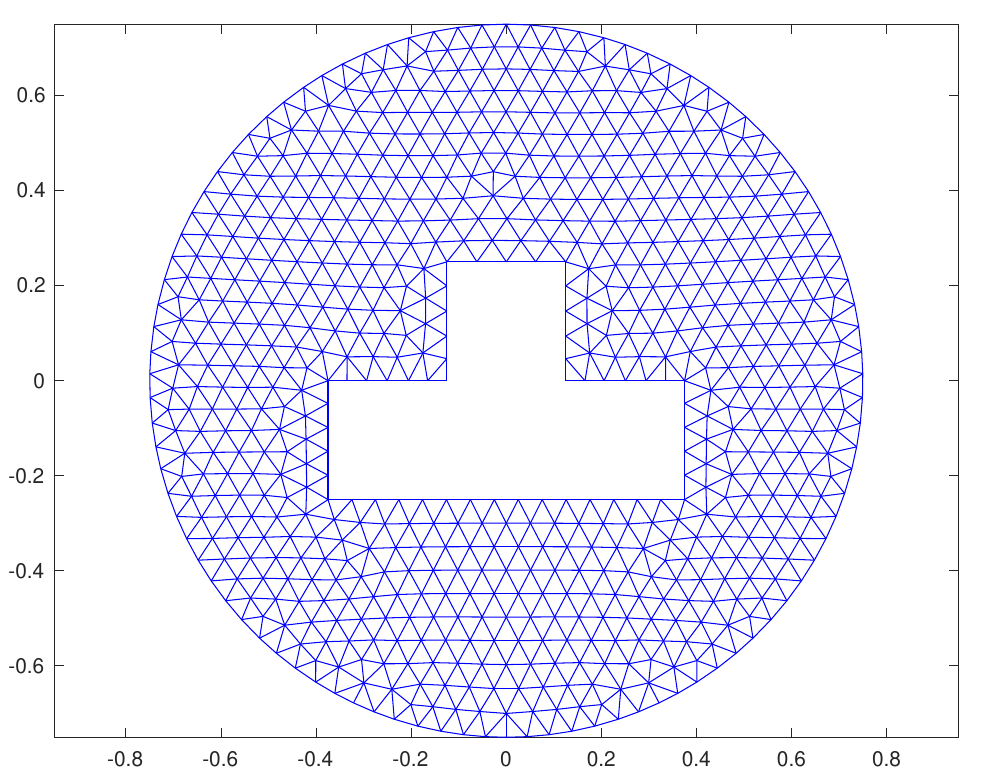}
}
\quad
\subfigure[$t=0.01$]{
\includegraphics[width=4.3cm]{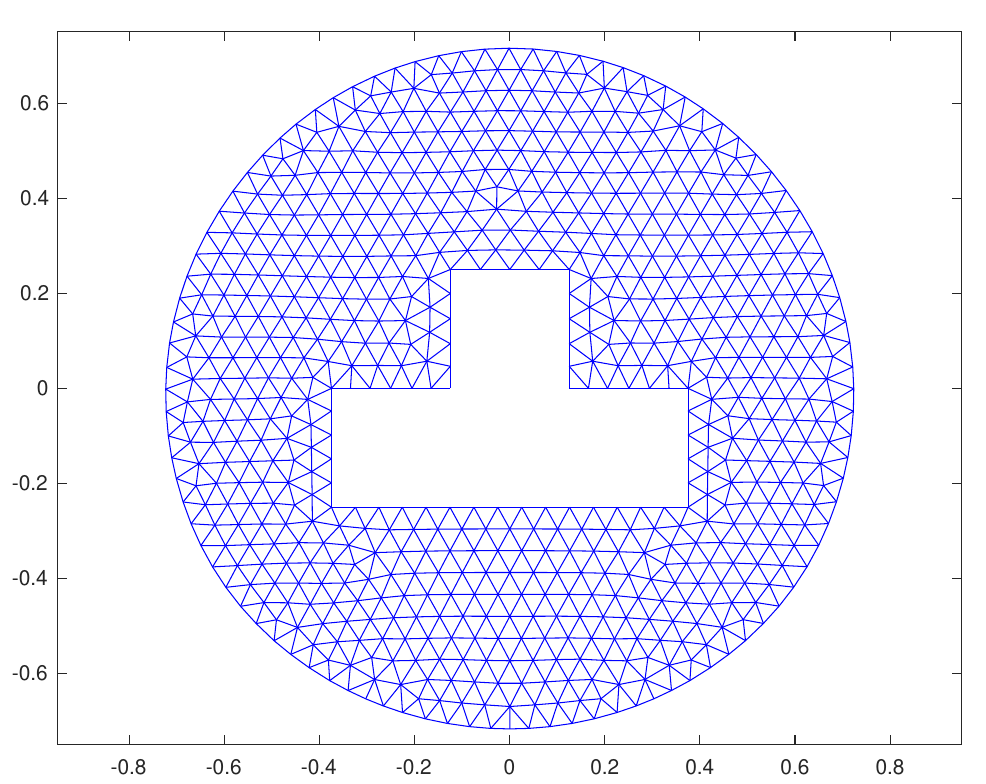}
}
\quad
\subfigure[$t=0.05$]{
\includegraphics[width=4.3cm]{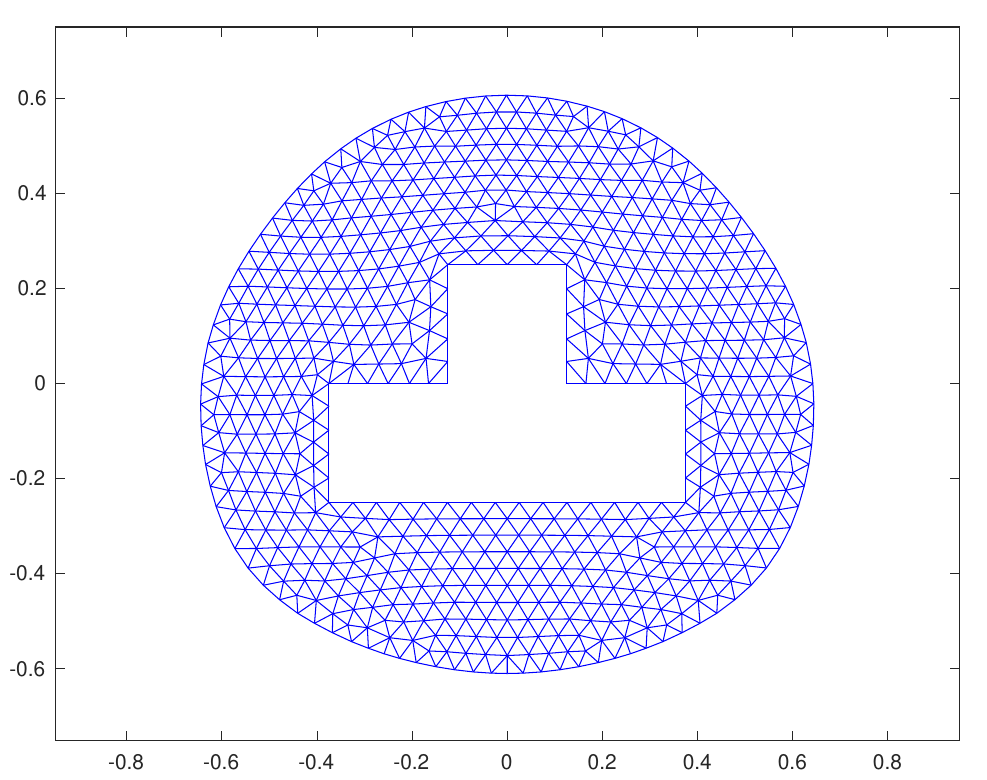}
}
\\
\subfigure[$t=0.15$]{
\includegraphics[width=4.4cm]{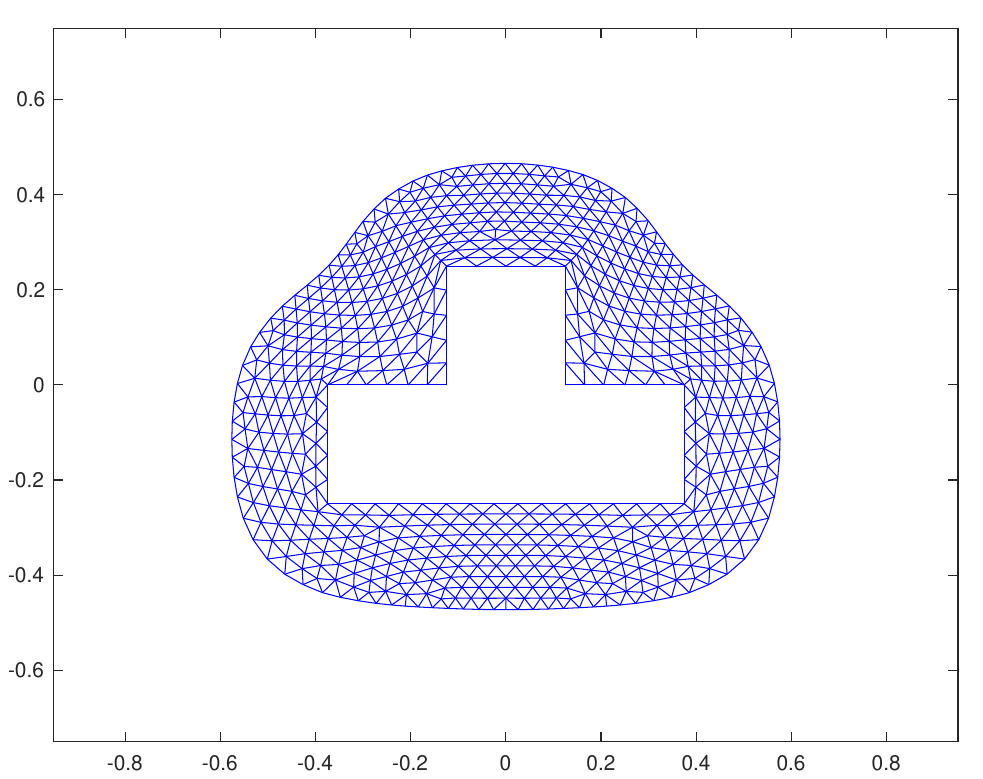}
}
\quad
\subfigure[$t=0.3$]{
\includegraphics[width=4.4cm]{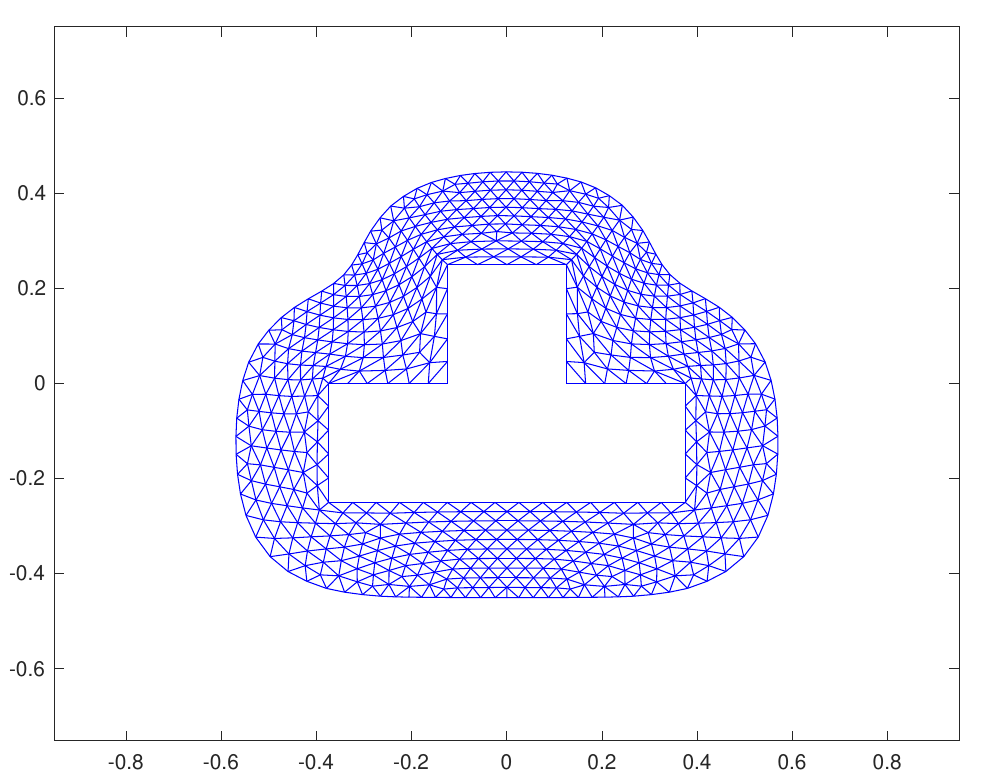}
}
\quad
\subfigure[$t=0.708$]{
\includegraphics[width=4.4cm]{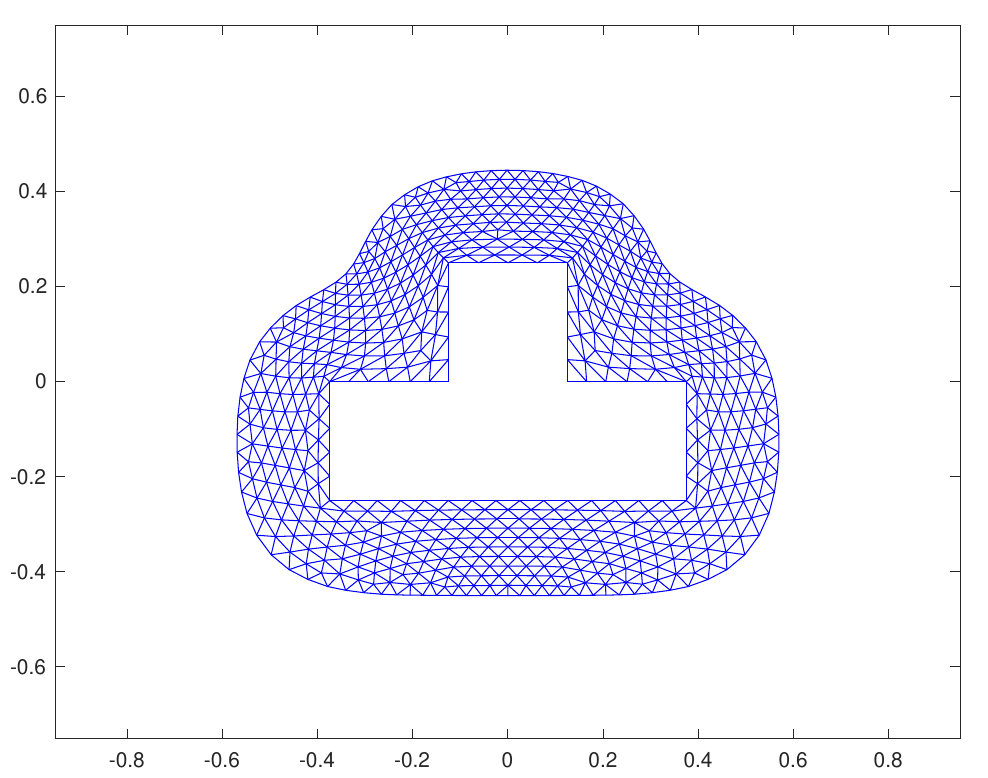}
}
\caption{Example~\ref{fbp-ex6} with $p = 5.0$. The mesh of $N = 1259$ is plotted at $t=0$, 0.05, 0.1, 0.15, 0.3, and 0.708 for $\lambda = 5$.}
\label{fig:fbp-ex6-2}
\end{figure}

\begin{figure}[htbp]
\centering
\subfigure[$p=1.5$]{
\includegraphics[width=4.3cm]{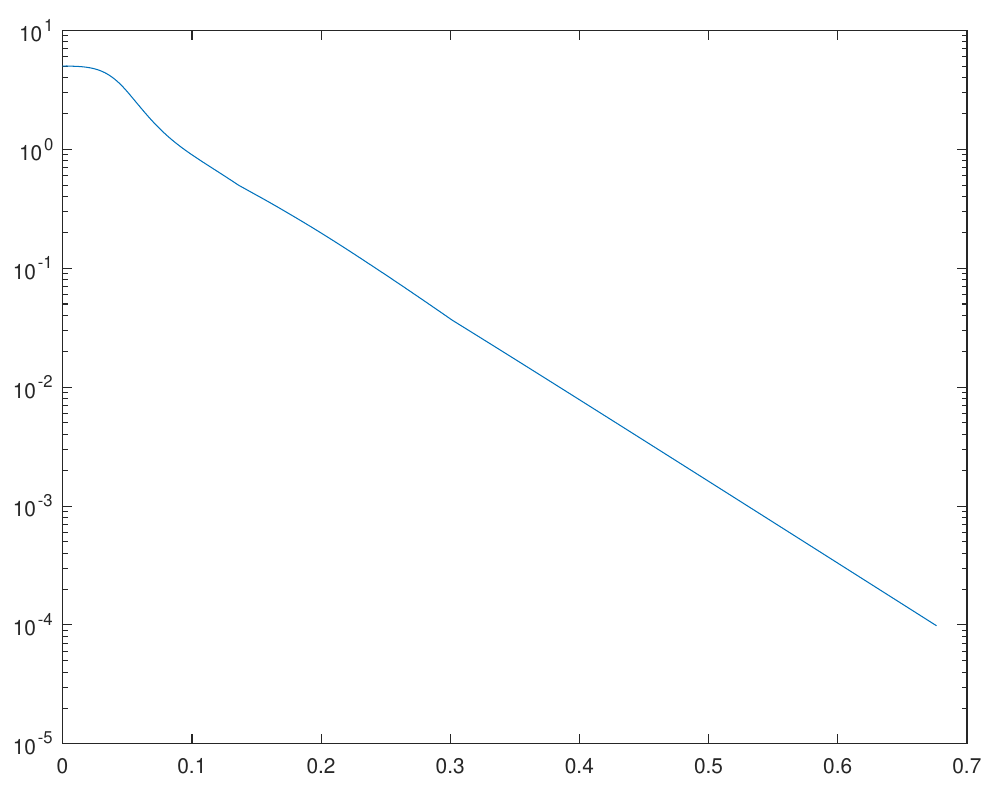}
}
\quad
\subfigure[$p=5.0$]{
\includegraphics[width=4.3cm]{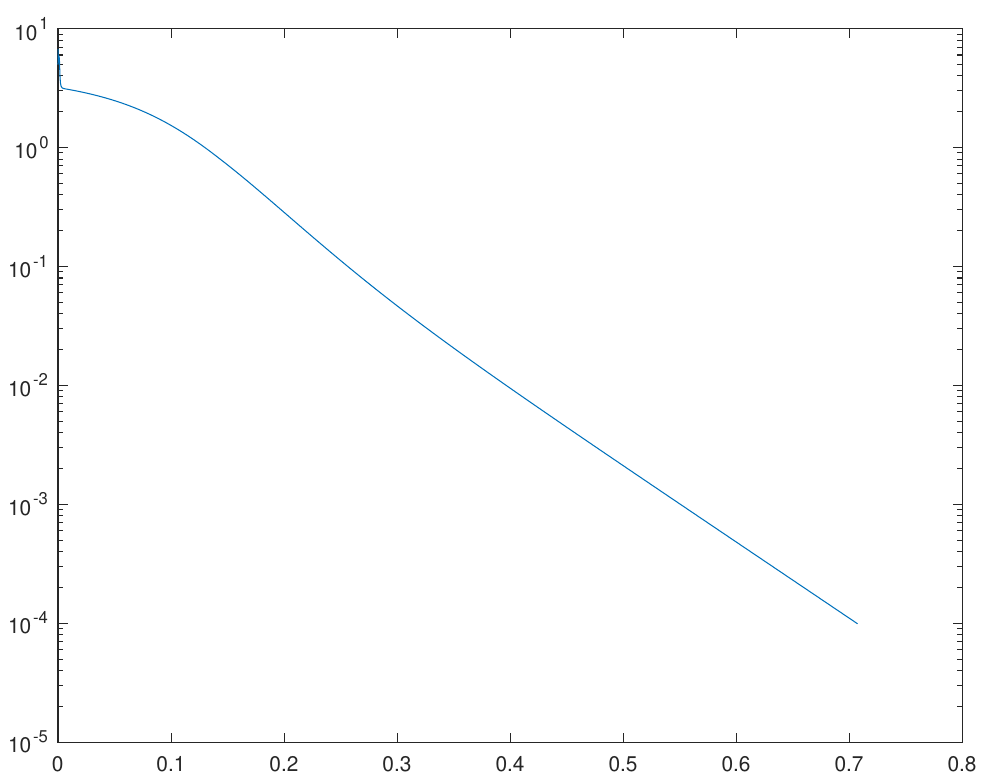}
}
\caption{Example~\ref{fbp-ex6}. The maximum boundary velocity is plotted as a function of time
for $\lambda = 5$, $N = 1259$, and two values of $p$.}
\label{fig:fbp-ex6-3}
\end{figure}

\begin{figure}[h!]
\centering
\subfigure[$p=1.5$]{
\includegraphics[width=4.3cm]{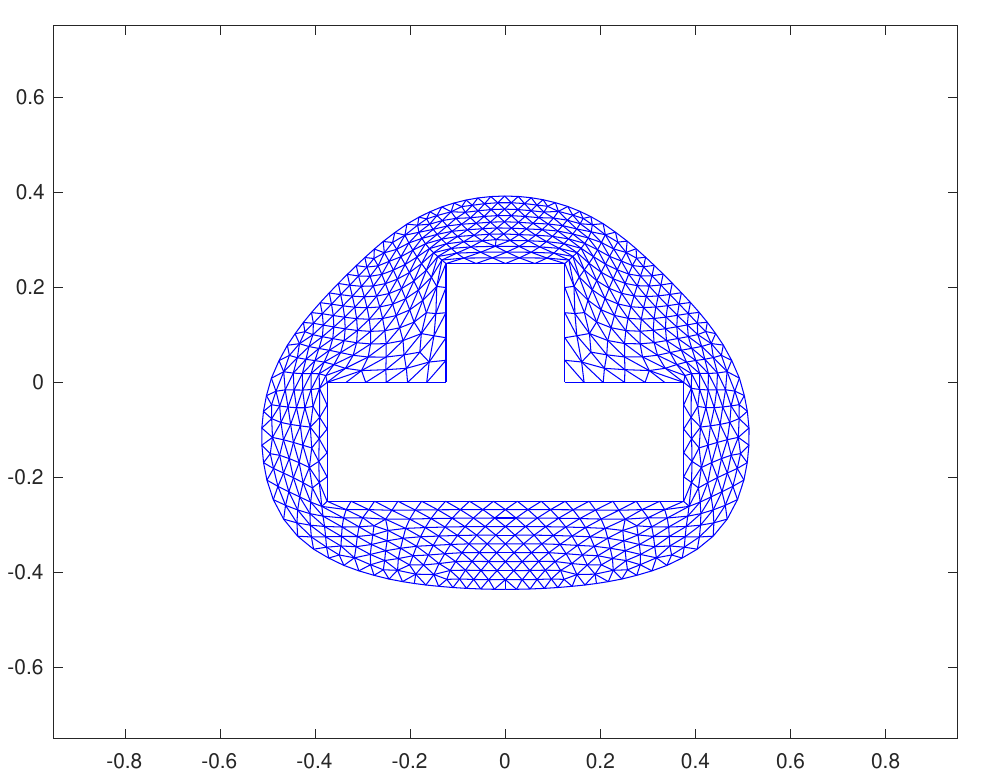}
}
\quad
\subfigure[$p=2.0$]{
\includegraphics[width=4.6cm]{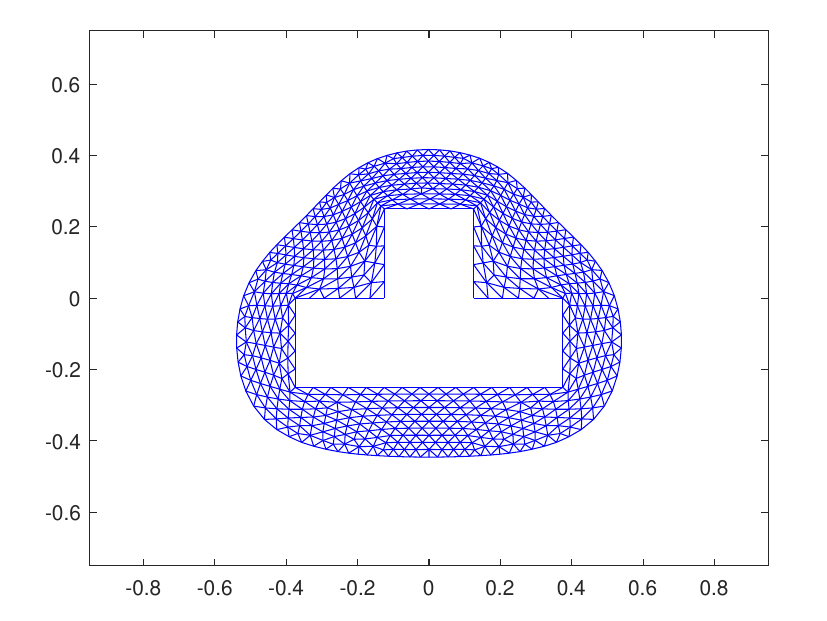}
}
\quad
\subfigure[$p = 5.0$]{
\includegraphics[width=4.3cm]{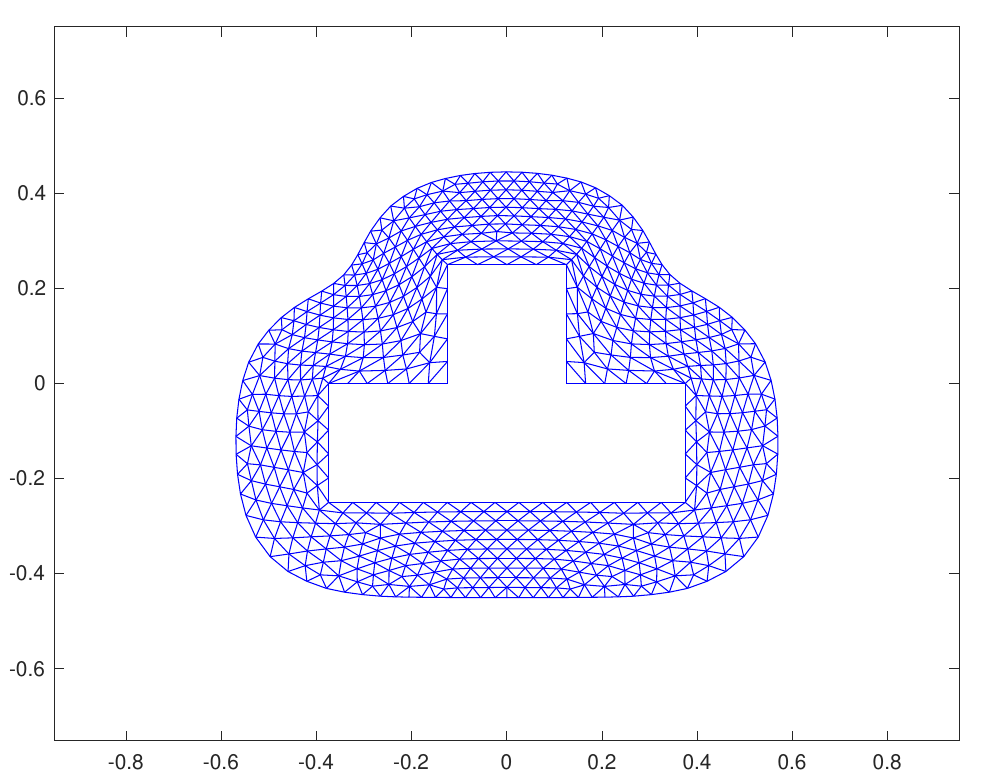}
}
\caption{Example~\ref{fbp-ex6} with $N = 1259$ and $\lambda = 5$. The mesh at $t = 0.3$ is compared for $p = 1.5$, $2.0$,
and $5.0$.}
\label{fig:fbp-ex6-4}
\end{figure}

\begin{exam}[\textbf{Exterior Bernoulli FBP with $p$-Laplacian}]
\label{fbp-ex6}
This example is the same as Example~\ref{fbp-ex2} except that the Laplace equation is replaced by the $p$-Laplace equation,
\begin{equation}
\nabla \cdot \big ( |\nabla u|^{p-2} \nabla u\big ) = 0,\quad \text{ in } \Omega
\label{p-Laplacian-1}
\end{equation}
where $p \in (1, \infty)$ is a parameter. The $p$-Laplacian is a power-law generalization of
various linear flow laws and is more realistic than the Laplacian (e.g., see Acker and Meyer \cite{Acker-1995}).
We take two values of $p$, $1.5$ and $5$, in our computation.
The meshes obtained with $p = 1.5$ and $p=5$ are shown in Figs.~\ref{fig:fbp-ex6-1} and \ref{fig:fbp-ex6-2}, respectively,
and the corresponding maximum boundary velocities are plotted in Fig.~\ref{fig:fbp-ex6-3}.
They confirm that the moving mesh FEM together with the pseudo-transient continuation works well for this nonlinear
example. Moreover, Fig.~\ref{fig:fbp-ex6-4} shows that the steady-state position of $\Gamma_2$ is more uniformly close to
$\Gamma_1$ for larger $p$.
\qed
\end{exam}



\begin{figure}[h!]
\centering
\subfigure[mesh at $t = 0$]{
\includegraphics[width=4.3cm]{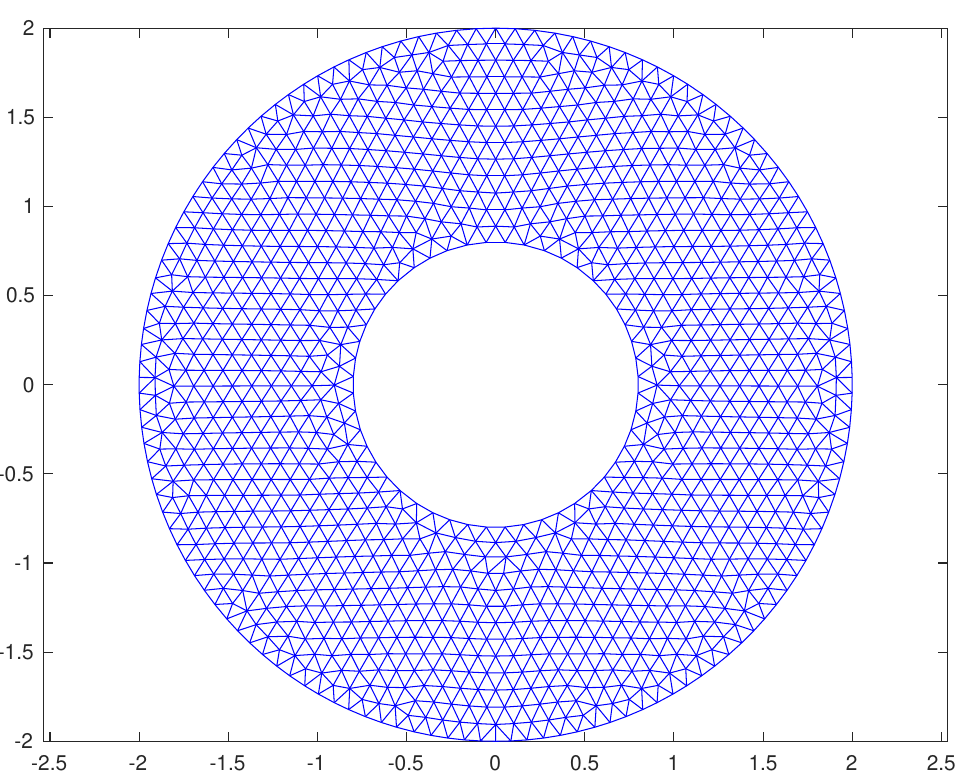}
}
\quad
\subfigure[mesh at $t = 2.686$]{
\includegraphics[width=4.3cm]{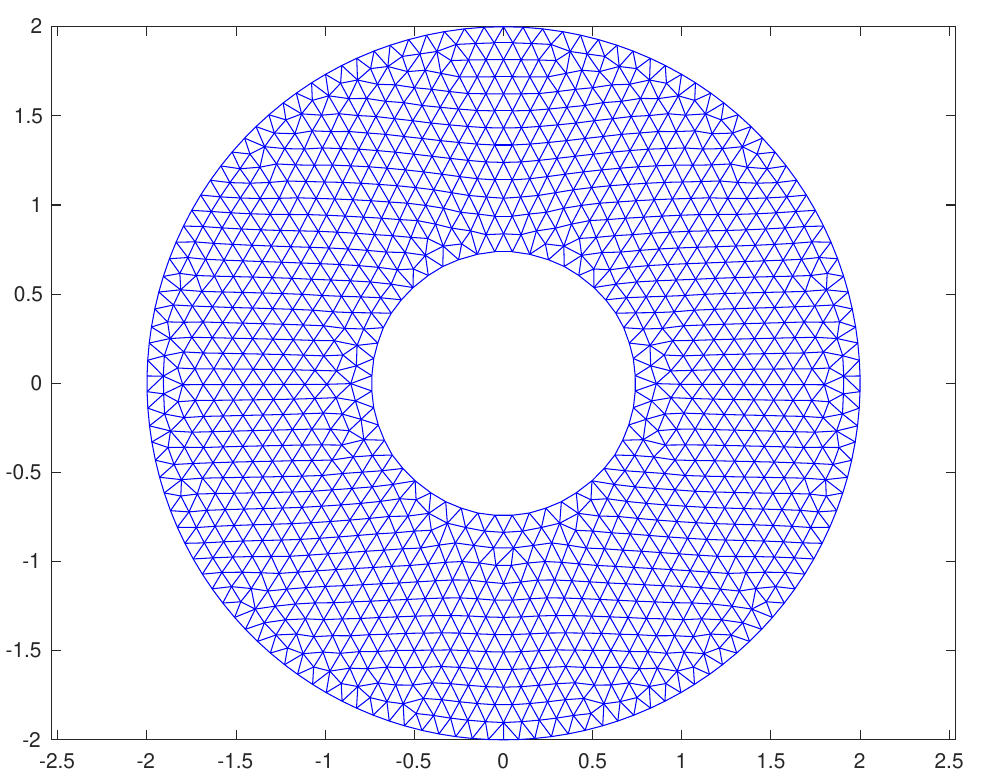}
}
\\
\subfigure[solution at $t = 2.686$]{
\includegraphics[width=4.3cm]{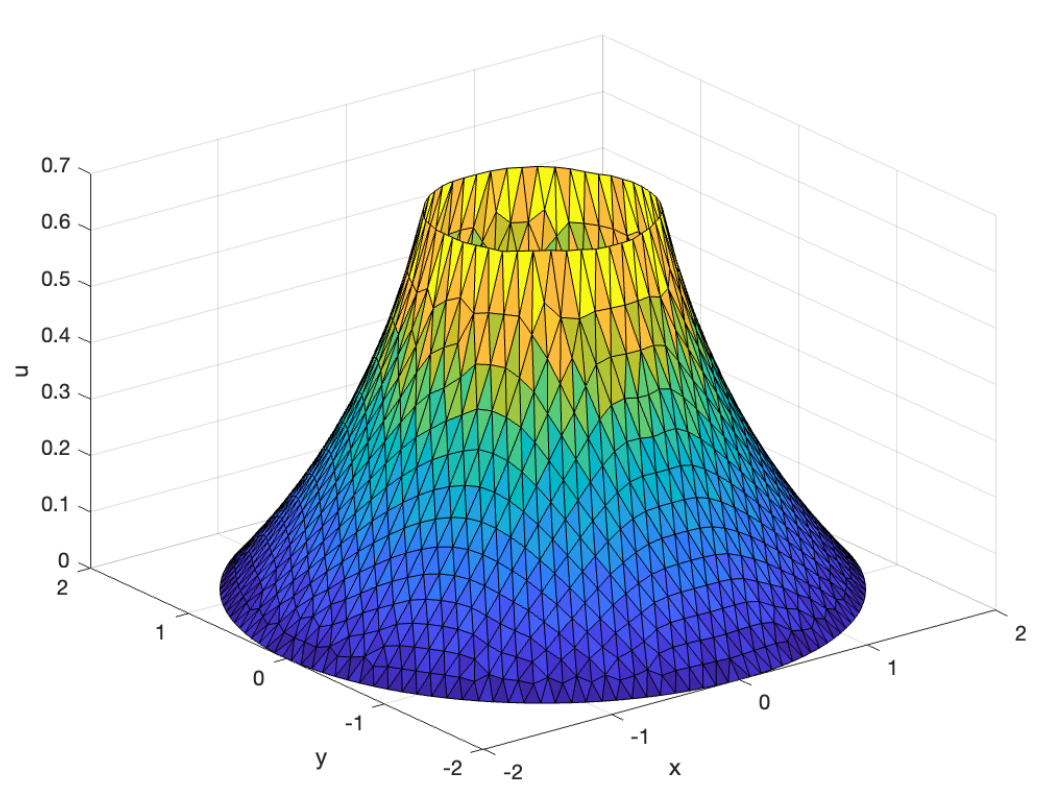}
}
\quad
\subfigure[max. boundary velocity vs $t$]{
\includegraphics[width=4.3cm]{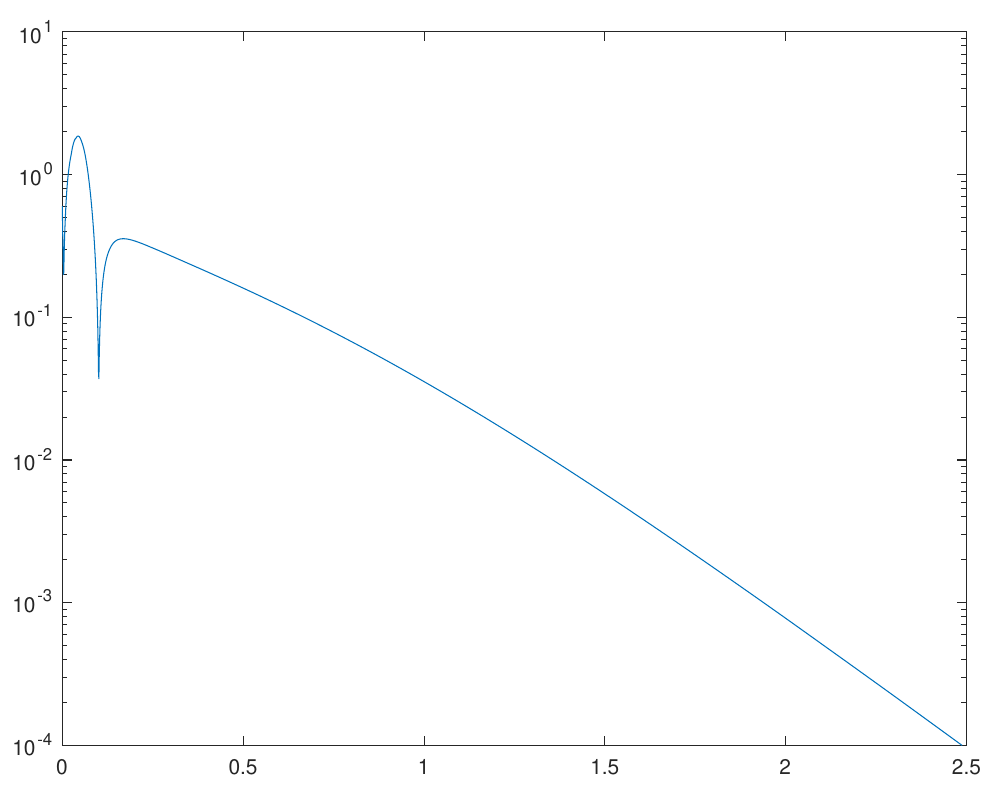}
}
\caption{Example~\ref{fbp-ex7}. The results are obtained with a mesh of $N = 2264$.}
\label{fig:fbp-ex7}
\end{figure}

\begin{exam}[\textbf{A nonlinear obstacle problem}]
\label{fbp-ex7}
Obstacle problems are a classical and important type of free boundary problem where the solution can be thought
as the equilibrium position of an elastic membrane that is constrained
to lie above a given obstacle $\psi = \psi(\V{x})$ while its boundary is held fixed (e.g., see Ros-Oton \cite{Ros-Oton-2018}).
We consider here a nonlinear obstacle problem
\begin{equation}
\label{fbp-ex7-1}
\min_{u} \int_{D} \sqrt{1 + |\nabla u|^2} d \V{x},\quad \text{subject to} \quad u \ge \psi \text{ in } D, \quad u = \psi \text{ on } \partial D
\end{equation}
where $D$ is the disk with radius $2$ and
\[
\psi = \begin{cases} \sqrt{1 - x^2 - y^2}, & \quad x^2 + y^2 \le 1 \\ 0, & \quad \text{otherwise}. \end{cases}
\]
This problem can be reformulated into a free boundary problem as
\begin{equation}
\label{fbp-ex7-2}
\begin{cases}
- \nabla \cdot \left (\frac{1}{\sqrt{1 + |\nabla u|^2}} \nabla u \right ) = 0, & \quad \text{ in } \Omega
\\
u = \psi, & \quad \text{ on } \Gamma_1 = \partial D
\\
u = \psi, & \quad \text{ on } \Gamma_2
\\
\frac{\partial u}{\partial n} = \frac{\partial \psi}{\partial n}, & \quad \text{ on } \Gamma_2
\end{cases}
\end{equation}
where $\Gamma_2$ is a closed curve inside $D$, $\Omega = D \setminus \overline{E}$, and $E$ is the domain enclosed
by $\Gamma_2$.  The Neumann boundary condition on $\Gamma_2$ is mathematically equivalent to a Bernoulli condition
$| \nabla (u-\psi) | = 0$. Moreover, the corresponding MBP in the pseudo-transient continuation is given by
\begin{equation}
\label{fbp-ex7-3}
\begin{cases}
\frac{\partial u}{\partial t} = \nabla \cdot \left (\frac{1}{\sqrt{1 + |\nabla u|^2}} \nabla u \right ), & \quad \text{ in } \Omega
\\
u = \psi, & \quad \text{ on } \Gamma_1 = \partial D
\\
u = \psi, & \quad \text{ on } \Gamma_2
\\
\dot{\Gamma} = - \frac{\partial u}{\partial n} + \frac{\partial \psi}{\partial n}, & \quad \text{ on } \Gamma_2 .
\end{cases}
\end{equation}
The initial condition for $u$ is taken as $u(\V{x},0) = 0$ and the initial position of $\Gamma_2$ is chosen as the circle with radius 0.8.
The mesh, solution, and maximum boundary velocity obtained with a mesh of $N = 2264$ are plotted in Fig.~\ref{fig:fbp-ex7}.
The results demonstrate that the moving mesh FEM and the pseudo-transient continuation can be used to obtain
the solution of the nonlinear obstacle problem (\ref{fbp-ex7-1}) as the steady-state solution of (\ref{fbp-ex7-3}).
\qed
\end{exam}
\pagebreak
\section{Conclusions and comments}
\label{SEC:conclusions}

We have studied a moving mesh finite element method for the numerical solution of Bernoulli FBPs.
The method is based on the pseudo-transient continuation with which an MBP is constructed
and its steady-state solution is taken as the solution of the underlying Bernoulli FBP.
The MBP is solved in a split manner at each time step: the moving boundary is updated with the Euler scheme,
the interior mesh points are moved using the MMPDE moving mesh method, and the corresponding
initial-boundary value problem is solved using the linear FEM.
The overall procedure is listed in Algorithm~\ref{MMFEM-euler}.
The method can take full advantages of both the pseudo-transient continuation and the MMPDE method.
Particularly, it is able to move the mesh, free of tangling, to fit the varying domain for a variety of geometries,
no matter if they are convex or concave. Moreover, it is convergent towards steady state for a broad class
of FBPs and initial guesses of the free boundary.

Numerical examples for Bernoulli FBPs with constant and non-constant Bernoulli conditions
and nonlinear FBPs have been presented. Numerical results have shown that the method
is second-order in space when the gradient of the solution at boundary vertices that is needed in free boundary update
is recovered with quadratic least squares fitting. Moreover, they have also shown that the method works well
for both exterior and interior Bernoulli FBPs with complex geometries and nonlinear FBPs.

Finally, we comment that while it is generally more robust than Newton's method, the pseudo-transient continuation is typically
slower than the latter (in terms of convergence towards steady state). Unfortunately, the moving mesh method studied in this work
also inherits this drawback from the pseudo-transient continuation. It is interesting to see how the method can be sped up.
One idea is to use a Davidenko-like equation or a preconditioner (e.g. see Kramer \cite{Kramer-1992}) when constructing
the moving boundary problem in the pseudo-transient continuation.
A main challenge on this is how to speed up the movement of the free boundary while avoiding mesh tangling.
Another issue is how to compute hyperbolic solutions \cite{Henrot-2021}.
As suggested by the formal analysis in Section~\ref{SEC:PTC}
or Fig.~\ref{fig:boundary-movement}, it seems that
the pseudo-transient continuation and thus the moving mesh FEM studied in this work can be used only for elliptic solutions.
It is interesting to see if a method based on the pseudo-transient continuation can be designed for computing hyperbolic solutions.

\vspace{20pt}

\section*{Acknowledgments}
  J. Shen was supported in part by the National Natural Science Foundation of China through grant [12101509] and W. Huang was supported in part by the University of Kansas General Research Fund FY23.


\end{document}